\documentclass[reqno,12pt]{amsart}
\usepackage[all,line,poly,rotate,ps,dvips]{xy}
\SelectTips{cm}{}

\usepackage{amssymb}
\usepackage{enumerate}

\numberwithin{equation}{section}

\newtheorem{theorem}[equation]{Theorem}

\newtheorem{corollary}[equation]{Corollary}
\newtheorem{claim}[equation]{Claim}
\newtheorem{lemma}[equation]{Lemma}
\newtheorem{proposition}[equation]{Proposition}

\newtheorem{thm}[equation]{Theorem}

\newtheorem{lem}[equation]{Lemma}
\newtheorem{prop}[equation]{Proposition}

\theoremstyle{definition}

\newtheorem{definition}[equation]{Definition}
\newtheorem{remark}[equation]{Remark}
\newtheorem{example}[equation]{Example}
\newtheorem{notation}[equation]{Notation}
\newtheorem{constr}[equation]{Construction}

\theoremstyle{remark}

\newcommand{\C}{\mathbb{C}}
\newcommand{\G}{\mathbb{G}}
\newcommand{\pp}{\mathbb{P}}
\newcommand{\Pp}{\mathbb{P}}
\newcommand{\F}{\mathbb{F}}
\newcommand{\Ff}{\mathcal{F}}
\newcommand{\CC}{\mathcal{C}}
\newcommand{\s}{\mathcal{S}}
\newcommand{\HH}{\mathcal{H}}

\newcommand{\II}{\mathcal{I}}
\newcommand{\Ii}{\mathcal{I}}

\newcommand{\N}{\mathcal{N}}
\newcommand{\T}{\mathcal{T}}
\newcommand{\Oc}{\mathcal{O}}
\newcommand{\OO}{\mathcal{O}}
\newcommand{\D}{\Delta}
\newcommand{\X}{\mathcal{X}}
\DeclareMathOperator{\rk}{rank}
\DeclareMathOperator{\coker}{coker}
\DeclareMathOperator{\Sing}{Sing}

\newcommand{\st}{\scriptstyle}
\newcommand{\dt}{\displaystyle}
\begin{document}
\title[Degenerations of scrolls to unions of planes]%
{Degenerations of scrolls to unions of planes}%
\author{A. Calabri, C. Ciliberto, F. Flamini, R. Miranda}

\email{calabri@dm.unibo.it}
\curraddr{Dipartimento di Matematica, Universit\`a degli Studi di Bologna\\
Piazza di Porta San Donato, 5 - 40126 Bologna \\Italy}

\email{cilibert@mat.uniroma2.it} \curraddr{Dipartimento di
Matematica, Universit\`a degli Studi di Roma Tor Vergata\\ Via
della Ricerca Scientifica - 00133 Roma \\Italy}

\email{flamini@mat.uniroma2.it} \curraddr{Dipartimento di
Matematica, Universit\`a degli Studi di Roma Tor Vergata\\ Via
della Ricerca Scientifica - 00133 Roma \\Italy}

\email{miranda@math.colostate.edu} \curraddr{Department of
Mathematics, 101 Weber Building, Colorado State University
\\ Fort Collins, CO 80523-1874 \\USA}

\thanks{
\\
The first three authors are members of G.N.S.A.G.A.\ at
I.N.d.A.M.\ ``Francesco Severi''.}

\begin{abstract} In this paper we study degenerations of scrolls to union of planes, a problem
already considered by G. Zappa in \cite{Zap1} and \cite{Zap2}. We prove, using techniques different from
the ones of Zappa, a degeneration result to union
of planes with the mildest possible singularities, for linearly
normal scrolls of genus $g$ and of degree $d \geq 2g+4$ in $\Pp^{d-2g+1}$. We also study properties of components
of the Hilbert scheme parametrizing scrolls. Finally we review Zappa's original approach.
\end{abstract}

\maketitle


\noindent
\centerline{{\em Dedicated to Professor G. Zappa on his 90th birthday}}

\section{Introduction}\label{Intro}

In this paper we deal with the problem,
originally studied by Guido Zappa in \cite{Zap1, Zap2},
concerning the embedded degenerations of two-dimensional scrolls,
to union of planes with the simplest possible singularities.

In \cite{CCFMFano} and  \cite{CCFMk2},
we have studied the properties of the so-called \emph{Zappatic surfaces}, i.e.
reduced, connected, projective
surfaces which are unions of smooth surfaces with global normal crossings
except at singular points, which are
locally analytically isomorphic to the vertex of a cone over a union of lines
whose dual graph is either a chain of length $n$, or a fork with $n-1$ teeth,
or a cycle of order $n$, and with maximal embedding dimension.
These singular points are respectively called  \emph{(good) Zappatic singularities} of type $R_n$, $S_n$ and $E_n$
(cf.\ Definition \ref{def:zapp} below).
A Zappatic surface is said to be \emph{planar} if it is embedded in a projective space and
all its irreducible components are planes.

An interesting problem is to find degenerations of surfaces to Zappatic surfaces
with Zappatic singularities as simple as possible. This problem has been partly considered in
\cite{CCFMk2}; e.g. in Corollary 8.10,
it has been shown that, if $X$ is a Zappatic surface which is the flat limit of a smooth
scroll of sectional genus $g \geq 2$, then the Zappatic singularities of $X$ cannot be too simple,
in particular $X$ has to have some point of type $R_i$ or $S_i$, with $i \geq 4$, or of type $E_j$,
with $j \geq 6$.

The main results in \cite{Zap1} can be stated in the following way:

\begin{theorem}\label{thm:zappa} (cf.\ \S 12 in \cite{Zap1})
Let $F$ be a scroll of sectional genus $g$, degree $d\geq 3g+2$, whose general hyperplane section
is a general curve of genus $g$. Then $F$ is birationally equivalent to a scroll in $\Pp^r$, for some $r \geq 3$,
which degenerates to a planar Zappatic surface with only points of
type $R_3$ and $S_4$ as Zappatic singularities.
\end{theorem}

Zappa's arguments rely on a rather intricate analysis concerning
degenerations of hyperplane sections of the scroll and, accordingly, of the branch curve of a general
projection of the scroll to a plane.

We have not been able to check all the details of this very clever argument. However,
we have been able to prove a slightly more general result using  some basic
smoothing technique (cf.\ \cite{CLM}).

Our main result is the following (cf.\ Proposition \ref{prop:rat1}, Constructions \ref{ex:3constr1s},
\ref{ex:3constr}, Remarks \ref{rem:ex6}, \ref{rem:hilbconto} and Theorems \ref{thm:deg3}, \ref{thm:contoflam} later on):

\begin{theorem}\label{thm:intro} Let $g \geq 0$ and either $d \geq 2$, if $g = 0$,  or
$d \geq  5$, if $g = 1$, or $d \geq 2g+4$, if $g \geq 2$. Then
there exists a unique irreducible component $\HH_{d,g}$ of the
Hilbert scheme of scrolls of degree $d$ and sectional genus $g$ in
$\Pp^{d-2g+1}$, such that the general point of $\HH_{d,g}$
represents a smooth scroll $S$ which is linearly normal and
moreover with $H^1(S,\Oc_S(1))=0$.

Furthermore,
\begin{itemize}
\item[(i)] $\HH_{d,g}$ is generically reduced and $dim (\HH_{d,g}) = (d-2g+2)^2+7(g-1)$,
\item[(ii)] $\HH_{d,g}$ contains the Hilbert point of
a planar Zappatic surface having only $d-2g+2$ points of type $R_3$ and
$2g-2$ points of type $S_4$ as Zappatic singularities,
\item[(iii)] $\HH_{d,g}$ dominates the moduli
space ${\mathcal M}_g$ of smooth curves of genus $g$.
\end{itemize}
\end{theorem}

We also construct examples of scrolls $S$ with same numerical invariants,
which are not linearly normal in $\Pp^{d-2g+1}$, as well as examples of components of the Hilbert scheme
of scrolls with same invariants, different from $\HH_{d,g}$ and with general moduli
(cf.\ Examples \ref{ex:contoflam2} and \ref{ex:ciro}).

We shortly describe the contents of the paper. In \S\ \ref{S:1bis}
we recall standard definitions and properties of Zappatic
surfaces. In \S\ \ref{S:ratscrolls} we focus on some degenerations
of products of curves to planar Zappatic surfaces and we prove
some results which go back to \cite{Zap2}. In particular, we
consider Zappatic degerations of rational and elliptic normal
scrolls and of abelian surfaces.

In \S\ \ref{S:4} we prove the greatest part of Theorem
\ref{thm:intro}. First, we construct, with an inductive argument,
planar Zappatic surfaces which have the same numerical invariants
of scrolls of degree $d$ and genus $g$ in $\Pp^{d-2g+1}$ and
having only  $d-2g+2$ points of type $R_3$ and $2g-2$ points of
type $S_4$ as Zappatic singularities. Then we prove that these
Zappatic surfaces can be smoothed to smooth scrolls which fill up
the component $\HH_{d,g}$ and we compute the cohomology of the
hyperplane bundle and of the normal bundle. These computations
imply that $\HH_{d,g}$ is generically smooth, of the right
dimension and its general point represents a linearly normal
scroll.

Section \ref{S:compo} is devoted to study some properties of
components of the  Hilbert scheme of scrolls. In particular, we
show that the component $\HH_{d,g}$ is the unique component of the
Hilbert scheme of scrolls of degree $d$ and sectional genus $g$
whose general point $[S]$ is linearly normal in $\Pp^{d-2g+1}$ and
moreover with $H^1(S, \Oc_S(1)) = 0$. Furthermore, we give the
examples mentioned above (cf.\ Examples \ref{ex:contoflam2} and
\ref{ex:ciro}).

In the last section, \S\ \ref{S:Zappa}, we briefly explain Zappa's original approach in \cite{Zap1}.
Moreover, we
make some comments and give some improvements on some interesting results from \cite{Zap1}
concerning extendability of plane curves to scrolls which are not cones.

\section{Notation and preliminaries}\label{S:1bis}

In this paper we deal with projective varieties defined over
the complex field $\mathbb{C}$.

Let us recall the notions of Zappatic singularities, Zappatic
surfaces and their dual graphs. We refer the
reader for more details to our previous papers \cite{CCFMFano} and
\cite{CCFMk2}. One word of warning: what we call {\em good
Zappatic singularities} there, here we simply call {\em Zappatic
singularities}, because no other type of Zappatic singularity will
be considered in this paper.

\begin{figure}[ht]
\vspace{-4.5em}
\[
\begin{xy}
<3.5em,0em>:
0; {\xypolygon10{~:{(1.75,0):}~>{}}},
(0,.4)="a"*=0{\bullet};"7"**@{-}?(.6)*{}="b" ,
"a";"8"**@{-}?(.6)*{}="c" ,
"a";"9"**@{-}?(.6)*{}="d" ,
"a";"10"**@{-}?(.6)*{}="e" ,
"b";"c"**@{-}?*!DL{\st } ,
"c";"d"**@{-}?+<0pt,4pt>*{\st } ,
"d";"e"**@{-}?*!DR{\st } ,
"7";"8"**@{}?*{V_1} ,
"8";"9"**@{}?*{V_2} ,
"9";"10"**@{}?*{V_3} ,
"8"*!U{\st C_{12}} ,
"9"*!U{\st C_{23}} ,
\end{xy}
\begin{xy}
<3.5em,0em>:
0; {\xypolygon10{~:{(1.75,0):}~>{}}},
(0,.4)="a"*=0{\bullet};"7"**@{-}?(.6)="b" ,
"a";"8"**@{-}?(.6)="c" ,
"a";"9"**@{-}?(.6)="d" ,
"8";"9"**@{}?="m" ,
"a";"10"**@{-}?(.6)="e" ,
"a";"m"**@{-} ?!{"c";"d"}="n" ,
"b";"c"**@{-}?*!DL{\st } ,
"c";"n"**@{-}?+<0pt,4pt>*{\st } ,
"d";"e"**@{-}?*!DR{\st } ,
"7";"8"**@{}?*{V_1} ,
"8";"m"**@{}?+<3pt,10pt>*!D{V_2} ,
"9";"10"**@{}?*{V_3} ,
"8"*!U{\st C_{12}} ,
"9"+<2pt,0pt>*!U{\st C_{23}} ,
"m";"9"**@{}?(.8)-(0,0.4)="p" ,
"a";"p"**@{-} ?(.6)="q" ,
"m";"9"**@{}?(.35)+<0pt,1pt>*!D{\dt V_4} ,
"m"*!U{\st C_{24}} ,
"n";"q"**@{-}?-<0pt,5pt>*{\st } ,
"n";"d"**@{} ?!{"a";"p"}="r" ,
"n";"r"**@{.} ,
"r";"d"**@{-} ,
\end{xy}
\begin{xy}
<3.5em,0em>:
0; {\xypolygon10{~={18}~:{(1.75,0):}~>{}}},
(0,.4)="a"*=0{\bullet};"7"**@{-}?(.6)*{}="b" ,
"a";"8"**@{-}?(.6)*{}="c" ,
"a";"9"**@{-}?(.6)*{}="d" ,
"b";"c"**@{-}?+<-1pt,-1pt>*!U{\st } ,
"c";"d"**@{-}?+<1pt,-1pt>*!U{\st } ,
"b";"d"**@{.}?(.35)*!D{\st } ,
"7";"8"**@{}?*{V_1} ,
"8";"9"**@{}?*{V_2} ,
"7";"9"**@{}?*+!DR{V_3} ,
"7"*!U{\st C_{13}} ,
"8"*!U{\st C_{12}} ,
"9"*!U{\st C_{23}} ,
\end{xy}
\]
\vspace{-6mm}
\caption{Examples: a $R_3$-point, a $S_4$-point and an $E_3$-point.}
\label{fig:R3E3S4}
\end{figure}
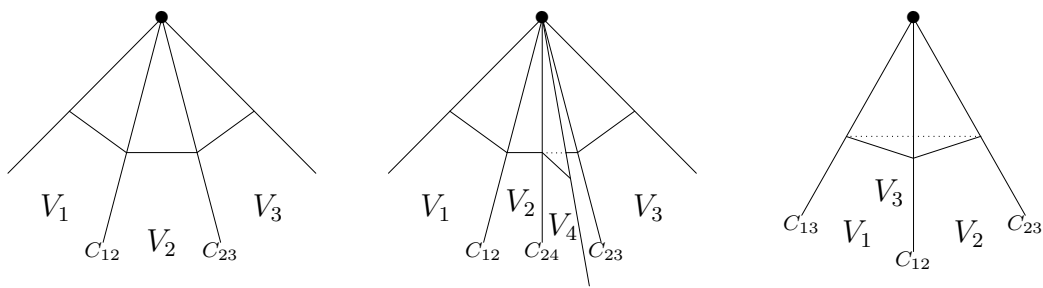

\begin{definition} \label{def:zapp}
Let us denote by $R_n$ [resp.\ $S_n$, $E_n$]
a graph which is a chain [resp.\ a fork, a cycle] with $n$ vertices, $n\geq3$,
cf.\ Figure \ref{fig:zapsings}.
Let $C_{R_n}$ [resp.\ $C_{S_n}$, $C_{E_n}$]
be a connected, projectively normal curve of degree $n$ in $\Pp^{n}$
[resp.\ in $\Pp^{n}$, in $\Pp^{n-1}$],
which is a {\em stick curve}, i.e.\ a reduced,
union of lines with only double points, whose dual graph is $R_n$
[resp.\ $S_n$, $E_n$].

\begin{figure}[ht]
\[
\raisebox{10pt}{$ %
\begin{xy}
0; <2em,0em>: 
(0,0)*=0{\bullet};(1,0)*=0{\bullet}**@{-};(2,0)*=0{\bullet}**@{-};(3,0)*=0{\bullet}**@{-};
(4,0)*=0{\bullet}**@{-};(5,0)*=0{\bullet}**@{-};(6,0)*=0{\bullet}**@{-}
,
\end{xy}
$}
\qquad
\qquad
\begin{xy}
0; <2em,0em>: 
(0,1.5)*=0{\bullet}="u" , (-2.5,0)*=0{\bullet};"u"**@{-} ,
(-1.5,0)*=0{\bullet};"u"**@{-} , (-0.5,0)*=0{\bullet};"u"**@{-} ,
(2.5,0)*=0{\bullet};"u"**@{-} , (1.5,0)*=0{\bullet};"u"**@{-} ,
(0.5,0)*=0{\bullet};"u"**@{-} ,
\end{xy}
\qquad
\qquad
\raisebox{15pt}{$ %
\begin{xy}
0; <2.5em,0em>: 
{\xypolygon7{\bullet}}
\end{xy}
$}
\]
\vspace{-1em}
\caption{A chain $R_n$, a fork $S_n$ with $n-1$
teeth, a cycle $E_n$.}\label{fig:zapsings}
\end{figure}
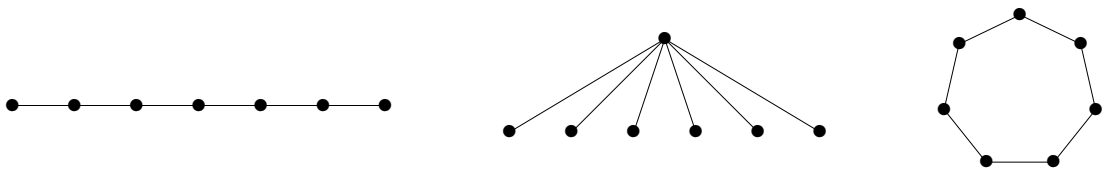

We say that a point $x$ of a projective surface $X$
is a \emph{point of type $R_n$} [resp.\ \emph{$S_n$, $E_n$}]
if $(X,x)$ is locally analytically isomorphic to a pair $(Y,y)$
where $Y$ is the cone over a curve $C_{R_n}$ [resp.\ $C_{S_n}$, $C_{E_n}$],
$n\geq 3$,
and $y$ is the vertex of the cone (cf.\ Figure \ref{fig:R3E3S4}).
We say that $R_n$-, $S_n$-, $E_n$-points are \emph{Zappatic singularities}.
\end{definition}

In this paper we will deal mainly with points of type $R_3$ and $S_4$.
We will use the following:

\begin{notation}\label{nota:central}
If $x$ is a point of type $R_3$ [of type $S_4$, resp.] of a projective surface $X$, we say that the component $V_2$ of $X$
as in picture on the left [in the middle, resp.] in Figure \ref{fig:R3E3S4} is the \emph{central} component of $X$
passing through $x$.
\end{notation}

\begin{definition}\label{def:zappsurf}
A projective surface $X=\bigcup_{i=1}^v V_i$ is called a \emph{Zappatic
surface} if $X$ is connected, reduced, all its irreducible components $V_1,
\ldots, V_v$ are smooth and:
\begin{itemize}
\item the singularities in codimension one of $X$ are at most
double curves which are smooth and irreducible along which two surfaces
meet transversally;
\item the further singularities of $X$ are Zappatic singularities.
\end{itemize}
We set $C_{ij}=V_i \cap V_j$ if $V_i$ and $V_j$ meet along a
curve, we set $C_{ij}=\emptyset$ otherwise. We set $C_i=V_i\cap
\overline {X-V_i}= \bigcup_{j=1}^ v C_{ij}$. We denote by $C=\Sing(X)$ the singular
locus of $X$, i.e.\ the curve $C= \bigcup_{1\leq i<j\leq v} C_{ij}$.

We denote by $f_n$ [resp.\ $r_n$, $s_n$] the number of point of type $E_n$ [resp.\ $R_n$, $S_n$] of $X$.
\end{definition}

\begin{remark}\label{rem:pg}
A Zappatic surface $X$ is Cohen-Macaulay. Moreover it has
global normal crossings except at the $R_n$- and $S_n$-points, for
$n\geq3$, and at the $E_m$-points, for $m\geq 4$.
\end{remark}

We associate to a Zappatic surface $X$ a \emph{dual graph} $G_X$ as follows.

\begin{definition} \label{def:dualgraph}
Let $X=\bigcup_{i=1}^v V_i$ be a Zappatic surface. The dual graph
$G_X$ of $X$ is given by:
\begin{itemize}

\item a vertex $v_i$ for each irreducible component $V_i$ of $X$;

\item an edge $l_{ij}$, joining the vertices $v_i$ and $v_j$, for
each irreducible component of the curve $C_{ij}=V_i\cap V_j$;

\item a \emph{$n$-face} $F_p$ for each point $p$ of $X$ of type $E_n$ for
some $n\geq 3$: the $n$ edges bounding the face $F_p$ are the $n$
irreducible components of the double curve $C$ of $X$ concurring at $p$;

\item an \emph{open $n$-face} for each point $p$ of $X$ of type $R_n$ for
some $n\geq 3$; it is bounded by $n-1$ edges,
corresponding to the $n-1$ irreducible components
of the double curve of $X$ concurring at $p$,
and by a \emph{dashed} edge, which we add in order to join the two extremal vertices;

\item a \emph{$n$-angle} for each $p$ of $X$ of type $S_n$,
spanned by the $n-1$ edges that are the $n-1$ irreducible components
of the double curves of $X$ concurring at $p$.
\end{itemize}
By abusing notation, we will denote by $G_X$ also the CW-complex
associated to the dual graph $G_X$ of $X$, formed by vertices, edges and $n$-faces.
\end{definition}

\begin{remark}[cf.\ \cite{CCFMFano}] \label{open3}
When we deal with the dual graph
of a \emph{planar} Zappatic surface $X=\bigcup_{i=1}^v V_i$, we will not indicate open $3$-faces
with a dashed edge. Indeed, the graph itself shows where open $3$-faces are located.
\end{remark}

Some invariants of a Zappatic surface $X$ have been computed
in \cite{CCFMFano} and in \cite{CCFMpg}, namely
the Euler-Poincar\'e characteristic $\chi(\Oc_X)$, the \emph{$\omega$-genus} $p_\omega(X)=h^0(X,\omega_X)$,
where $\omega_X$ is the dualizing sheaf of $X$, and,
when $X$ is embedded in a projective space $\pp^r$,
the \emph{sectional genus} $g(X)$,
i.e.\ the arithmetic genus of a general hyperplane section of $X$.
In particular, for a planar Zappatic surface
(for the general case, see \cite{CCFMFano, CCFMpg}) one has:

\begin{proposition}
Let $X=\bigcup_{i=1}^v V_i$ be a planar Zappatic surface of degree $v$ in $\pp^r$
and denote by $e$ the degree of $C=\Sing(X)$, i.e.\ the number of double lines of $X$.
Then:
\begin{align}
g(X) &= e-v+1,  \\
p_\omega(X)&= h^0(X,\omega_X) = h_2(G_X,\C), \\
\chi(\Oc_X) &= \chi(G_X) = v-e+\sum_{i\geq3} f_i.
\end{align}
\end{proposition}

In this paper, a Zappatic surface will always be considered
as the central fibre of an embedded degeneration, in the following sense.

\begin{definition}\label{def:degen}
Let $\D$ be the spectrum of a DVR (or equivalently the complex unit disk).
A \emph{degeneration} of surfaces parametrized by $\D$ is a proper
and flat morphism $\pi:\X \to \Delta$ such that
each fibre $\X_t = \pi^{-1}(t)$, $t \neq 0$ (where 0 is the closed point
of $\D$), is a smooth, irreducible,
projective surface.
A degeneration $\pi:\X\to\D$ is said to be \emph{embedded}
in $\Pp^r$ if $\X\subseteq \D\times\Pp^r$ and the following diagram commutes:
\[
\xymatrix@C=0mm@R=5mm{\X \ar[d]_\pi & \subseteq & \Delta \times \Pp^r
\ar[dll]^{\text{pr}_1} \\ \Delta}
\]
\end{definition}

The invariants of the Zappatic surface $X=\X_0$, which is the central fibre
of an embedded degeneration $\X\to\D$, determine the invariants
of the general fibre $\X_t$, $t\ne0$, as we proved in \cite{CCFMFano, CCFMk2, CCFMpg}.
Again, we recall these results only for planar Zappatic surfaces
and we refer to our previous papers for the general case.

\begin{theorem} \label{thm:invXt}
Let $\X\to\D$ be an embedded degeneration in $\pp^r$ such that
the central fibre $X=\X_0$ is a planar Zappatic surface.
Then, for any $0\ne t \in \D$:
\begin{equation}
g(\X_t) = g(X),
\qquad
p_g(\X_t) = p_\omega(X),
\qquad
\chi(\Oc_{\X_t})=\chi(\Oc_X).
\end{equation}
Moreover the self-intersection $K^2_{\X_t}$ of a canonical divisor of $\X_t$ is:
\begin{equation}
K^2_{\X_t} = 9v-10e+\sum_{n\geq3} 2nf_n + r_3 + k,
\end{equation}
where $k$ depends on the presence of points of type $R_m$ and $S_m$, $m\geq4$:
\[
\sum_{m\geq4} (m-2)(r_m+s_m) \leq k \leq
\sum_{m\geq4} (2m-5)r_m+ \binom{m-1}{2} s_m.
\]
\end{theorem}

Finally, let us recall the construction of rational normal scrolls.

\begin{definition} \label{def:scroll}
Fix two positive integers $a,b$ and set $r=a+b+1$.
In $\pp^r$ choose two disjoint linear spaces $\pp^a$ and $\pp^b$.
Let $C_a$ [resp.\ $C_b$] be a smooth, rational normal curve of degree $a$ in $\pp^a$ [resp.\ of degree $b$ in $\pp^b$]
and fix an isomorphism $\phi: C_a \to C_b$.
Then, the union in $\pp^r$ of all the lines $\overline{p,\phi(p)}$, $p\in C_a$, is
a smooth, rational, projectively normal surface which is called \emph{scroll of type $(a,b)$} and it is denoted
by $S_{a,b}$. Such a scroll is said to be \emph{balanced} if either $b=a$ or $b=a+1$.

Another way to define a scroll is as the embedding of a Hirzebruch surface $\F_n$, $n\geq0$,
which is the minimal ruled surface over $\pp^1$ with a section of self-intersection $(-n)$.
Setting $F$ the ruling of $\F_n$ and $C$ a section such that $C^2=n$,
the linear system $|C+aF|$ embeds $\F_n$ in $\pp^{n+2a+1}$ as a scroll of type $(a,a+n)$, cf.\ e.g.\ \cite{GH}.
In particular a balanced scroll in $\pp^r$, $r\geq3$, is the embedding either of $\F_0=\pp^1\times\pp^1$
or of $\F_1$ depending on whether $r$ is odd or even.
\end{definition}

In the next section we will see, in particular, degenerations of
rational scrolls to a planar Zappatic surface.
In the subsequent section we will deal with scrolls of higher genus.

\section{Degenerations of product of curves and of rational scrolls}\label{S:ratscrolls}

Zappa suggested in \cite{Zap2} an interesting method for degenerating products of curves,
which also gives a degeneration of rational and elliptic
scrolls to planar Zappatic surfaces with only $R_3$-points.

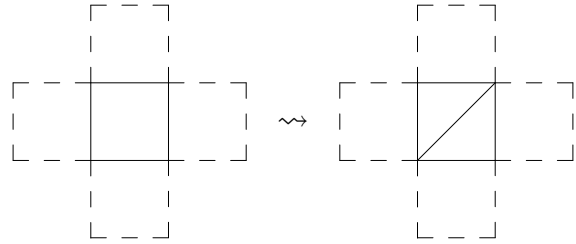
\begin{figure}[ht]
\[
\begin{xy}
0; <2.5em,0em>:
(1,1);(0,1)**@{--};(0,2)**@{--};(1,2)**@{--};(1,3)**@{--};(2,3)**@{--};
(2,2)**@{--};(3,2)**@{--};(3,1)**@{--};(2,1)**@{--};(2,0)**@{--};
(1,0)**@{--};(1,1)**@{--};
(2,1)**@{-};(2,2)**@{-};(1,2)**@{-};(1,1)**@{-} ,
\end{xy}
\quad
\raisebox{3.75em}{$\rightsquigarrow$}
\quad
\begin{xy}
0; <2.5em,0em>:
(1,1);(0,1)**@{--};(0,2)**@{--};(1,2)**@{--};(1,3)**@{--};(2,3)**@{--};
(2,2)**@{--};(3,2)**@{--};(3,1)**@{--};(2,1)**@{--};(2,0)**@{--};
(1,0)**@{--};(1,1)**@{--};
(2,1)**@{-};(2,2)**@{-};(1,2)**@{-};(1,1)**@{-} ;(2,2)**@{-} ,
\end{xy}
\]
\caption{A quadric degenerating to the union of two planes}\label{fig:quadrics}
\end{figure}

\begin{example}[Zappa]
Let $C\subset \Pp^{n-1}$ and $C'\subset \Pp^{m-1}$ be smooth curves.
If $C$ and $C'$ may degenerate to stick curves,
then the smooth surface
\[
S= C \times C'\subset \Pp^{n-1}\times \Pp^{m-1} \subset \Pp^{nm-1},
\]
embedded via the Segre map,
degenerates to a Zappatic surface $Y$ in $\Pp^{nm-1}$ whose irreducible components are quadrics
and whose double curves are lines.

If it is possible to further, independently, degenerate each quadric of $Y$ to the union of two planes,
then one gets a degeneration of $S=C \times C'$ to a planar Zappatic surface.
This certainly happens if each quadric of $Y$ meets the other quadrics of $Y$
along a union of at most four lines, at most two from each ruling
(see Figure \ref{fig:quadrics}).
\end{example}

Therefore $S =C\times C'$ can degenerate to a planar Zappatic surface
if $C$ and $C'$ are either rational or elliptic normal curves,
since they degenerate to stick curves $C_{R_n}$ and $C_{E_n}$,
respectively. We will now describe these degenerations.

\begin{figure}[ht]
\[
\begin{xy}
0; <2.5em,0em>:
(0,0);(6,0)**@{-};(6,1)**@{-};(0,1)**@{-};(0,0)**@{-} ,
(1,0);(1,1)**@{-} ,
(2,0);(2,1)**@{-} ,
(3,0);(3,1)**@{-} ,
(4,0);(4,1)**@{-} ,
(5,0);(5,1)**@{-} ,
(1,0);(1,1)**@{-} ,
(0,0);(1,1)**@{-} ,
(1,0);(2,1)**@{-} ,
(2,0);(3,1)**@{-} ,
(3,0);(4,1)**@{-} ,
(4,0);(5,1)**@{-} ,
(5,0);(6,1)**@{-} ,
\end{xy}
\qquad
\raisebox{1.25em}{$ %
\begin{xy}
0; <1.25em,0em>: 
(0,0)*=0{\bullet};(1,0)*=0{\bullet}**@{-};
(2,0)*=0{\bullet}**@{-};(3,0)*=0{\bullet}**@{-};
(4,0)*=0{\bullet}**@{-};(5,0)*=0{\bullet}**@{-};(6,0)*=0{\bullet}**@{-};
(7,0)*=0{\bullet}**@{-};(8,0)*=0{\bullet}**@{-};(9,0)*=0{\bullet}**@{-};
(10,0)*=0{\bullet}**@{-};(11,0)*=0{\bullet}**@{-}
,
\end{xy}
$}
\]
\caption{Planar Zappatic surface of degree $2n$ with a chain as dual graph}
\label{fig:2nplanes}
\end{figure}
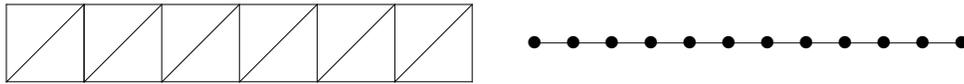

\begin{example}[Rational scrolls] \label{ex:zappa1}
Let $C$ be a smooth, rational normal curve of degree $n$ in $\pp^{n}$.
Since $C$ degenerates to a union of $n$ lines
whose dual graph is a chain,
the smooth rational normal scroll $S=C\times\pp^1\subset\pp^{2n+1}$
degenerates to a Zappatic surface $Y=\bigcup_{i=1}^n Y_i$ such that each $Y_i$ is a quadric,
$Y$ has no Zappatic singularity
and its dual graph $G_Y$ is a chain of length $n$,
see Figure \ref{fig:nquadrics}.

\begin{figure}[ht]
\[
\begin{xy}
0; <2.5em,0em>:
(0,0);(6,0)**@{-};(6,1)**@{-};(0,1)**@{-};(0,0)**@{-} ,
(1,0);(1,1)**@{-} ,
(2,0);(2,1)**@{-} ,
(3,0);(3,1)**@{-} ,
(4,0);(4,1)**@{-} ,
(5,0);(5,1)**@{-} ,
(1,0);(1,1)**@{-} ,
\end{xy}
\]
\caption{Chain of $n$ quadrics as in Example \ref{ex:zappa1}}
\label{fig:nquadrics}
\end{figure}
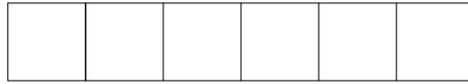

Each quadric $Y_i$ meets $\overline{Y\setminus Y_i}$ either along a line or along
two distinct lines of the same ruling.
Thus, as we noted before, the quadric $Y_i$ degenerates, in the $\pp^3$ spanned by $Y_i$,  to the union of two planes
meeting along a line $l_i$, leaving the other line(s) fixed.
Therefore, in $\pp^{2n+1}$, the scroll $S$ degenerates also to a planar Zappatic surface $X$ of degree $2n$.
The line $l_i$ can be chosen generally enough so that
$X$ has $2n-2$ points of type $R_3$ as Zappatic singularities, for each $i$,
i.e.\ its dual graph $G_X$ is a chain of length $2n$,
see Figure \ref{fig:2nplanes} (cf.\ Remark \ref{open3}).
\end{example}

\begin{example}[Elliptic scrolls] \label{ex:zappa2}
Let $C$ be a smooth, elliptic normal curve of degree $n$ in $\pp^{n-1}$.
Since $C$ degenerates to a union of $n$ lines
whose dual graph is a cycle,
the smooth elliptic normal scroll $S= C\times\pp^1\subset\pp^{2n-1}$
degenerates to a Zappatic surface $Y=\bigcup_{i=1}^n Y_i$, such that each $Y_i$ is quadric,
$Y$ has no Zappatic singularity
and its dual graph $G_Y$ is a cycle of length $n$,
see the picture on the left in Figure \ref{fig:e2n}.

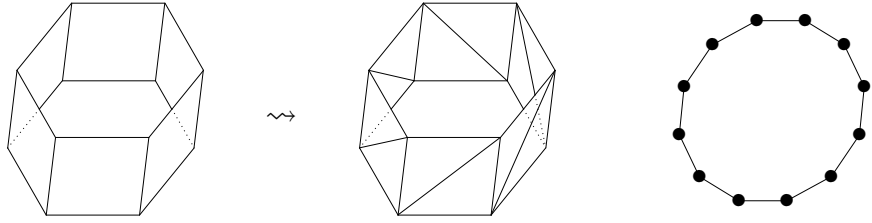
\begin{figure}[ht]
\[
\begin{xy}
0; <3em,0em>: <0.3em,2.5em>::
(0,1); {\xypolygon6"B"{~>{}}} ,
(0,0); {\xypolygon6"A"{~>{}}} ,
"A1";"A2"**@{-};"A3"**@{-};"A4"**@{-};"A5"**@{-};"A6"**@{-};"A1"**@{-} ,
"B1";"B6"**@{-};"B5"**@{-};"B4"**@{-} ,
"A1";"B1"**@{-} ,
"A2";"B2"**@{-} ,
"A3";"B3"**@{-} ,
"A4";"B4"**@{-} ,
"A5";"B5"**@{-} ,
"A6";"B6"**@{-} ,
"B1";"B2"**@{}?!{"A1";"A6"}="B12" ,
"B3";"B4"**@{}?!{"A4";"A5"}="B34" ,
"B1";"B12"**@{.};"B2"**@{-};"B3"**@{-};"B34"**@{-};"B4"**@{.} ,
\end{xy}
\qquad
\raisebox{1em}{$ \rightsquigarrow $}
\qquad
\begin{xy}
0; <3em,0em>: <0.3em,2.5em>::
(0,1); {\xypolygon6"B"{~>{}}} ,
(0,0); {\xypolygon6"A"{~>{}}} ,
"A1";"A2"**@{-};"A3"**@{-};"A4"**@{-};"A5"**@{-};"A6"**@{-};"A1"**@{-} ,
"B1";"B6"**@{-};"B5"**@{-};"B4"**@{-} ,
"A1";"B1"**@{-} ,
"A2";"B2"**@{-} ,
"A3";"B3"**@{-} ,
"A4";"B4"**@{-} ,
"A5";"B5"**@{-} ,
"A6";"B6"**@{-} ,
"B1";"B2"**@{}?!{"A1";"A6"}="B12" ,
"B3";"B4"**@{}?!{"A4";"A5"}="B34" ,
"B1";"B12"**@{.};"B2"**@{-};"B3"**@{-};"B34"**@{-};"B4"**@{.} ,
"B1";"A2"**@{}?!{"A1";"A6"}="AB1" ,
"B1";"AB1"**@{.};"A2"**@{-} ,
"B2";"A3"**@{-} ,
"B3";"A4"**@{-} ,
"B4";"A5"**@{-} ,
"B5";"A6"**@{-} ,
"B6";"A1"**@{-} ,
\end{xy}
\qquad
\qquad
\raisebox{-0.3em}{$ %
\begin{xy}
0; <3em,0em>: <0.3em,3em>::
(0,0.5) , {\xypolygon12{\bullet}}
\end{xy}
$}
\]
\caption{Cycle of $n$ quadrics and of $2n$ planes as in Example \ref{ex:zappa2}} \label{fig:e2n}
\end{figure}

Each quadric $Y_i$ meets $\overline{Y\setminus Y_i}$ along two distinct lines $r_i, r'_i$ of the same ruling.
Hence, in the $\pp^3$ spanned by $Y_i$, the quadric $Y_i$ degenerates to the union of two planes
meeting along a line $l_i$, leaving $r_i,r'_i$ fixed.
Choosing again a general $l_i$ for each $i$, it follows that in $\pp^{2n-1}$ the scroll $S$
degenerates to a planar Zappatic surface $X$ of degree $2n$ with $2n$ points of type $R_3$ as Zappatic singularities
and its dual graph $G_{X}$ is a cycle of length $2n$, see Figure \ref{fig:e2n}.
\end{example}

\begin{example} [Abelian surfaces]
Let $C \subset \Pp^{n-1}$ and $C'\subset \Pp^{m-1}$ be
smooth, elliptic normal curves of degree respectively $n$ and $m$.
Then $C$ and $C'$ degenerate to the stick curves $C_{E_n}$ and $C_{E_m}$ respectively,
hence the abelian surface $S = C \times C'\subset \Pp^{nm-1}$ degenerates to
a Zappatic surface which is a union of $mn$ quadrics with only $E_4$-points as Zappatic
singularities, cf.\ e.g.\ the picture on the left in Figure \ref{fig:quadrics2},
where the top edges have to be identified with the bottom ones,
similarly the left edges have to be identified with the right ones.
Thus the top quadrics meet the bottom quadrics and
the quadrics on the left meet the quadrics on the right.

\begin{figure}[ht]
\[
\begin{xy}
0; <2em,0em>: (0,0);(5,0)**@{-} , (0,4);(5,4)**@{-} ,
(0,0);(0,4)**@{-} , (5,0);(5,4)**@{-} , (0,1);(5,1)**@{-} ,
(0,2);(5,2)**@{-} , (0,3);(5,3)**@{-} , (1,0);(1,4)**@{-} ,
(2,0);(2,4)**@{-} , (3,0);(3,4)**@{-} , (4,0);(4,4)**@{-} ,
(0,0);(1,0)**@{}, (0,4);(1,4)**@{}, (1,0);(2,0)**@{},
(1,4);(2,4)**@{}, (2,0);(3,0)**@{}, (2,4);(3,4)**@{},
(3,0);(4,0)**@{}, (3,4);(4,4)**@{}, (4,0);(5,0)**@{},
(4,4);(5,4)**@{}, (0,0);(0,1)**@{}, (5,0);(5,1)**@{},
(0,1);(0,2)**@{}, (5,1);(5,2)**@{}, (0,2);(0,3)**@{},
(5,2);(5,3)**@{}, (0,3);(0,4)**@{}, (5,3);(5,4)**@{},
\end{xy}
\qquad
\raisebox{4em}{$ \rightsquigarrow $}
\qquad
\begin{xy}
0; <2em,0em>:
(4,0);(5,1)**@{-} ,
(3,0);(5,2)**@{-} ,
(2,0);(5,3)**@{-} ,
(1,0);(5,4)**@{-} ,
(0,0);(4,4)**@{-} ,
(0,1);(3,4)**@{-} ,
(0,2);(2,4)**@{-} ,
(0,3);(1,4)**@{-} ,
(0,0);(5,0)**@{-} ,
(0,4);(5,4)**@{-} ,
(0,0);(0,4)**@{-} ,
(5,0);(5,4)**@{-} ,
(0,1);(5,1)**@{-} ,
(0,2);(5,2)**@{-} ,
(0,3);(5,3)**@{-} ,
(1,0);(1,4)**@{-} ,
(2,0);(2,4)**@{-} ,
(3,0);(3,4)**@{-} ,
(4,0);(4,4)**@{-} ,
(0,0);(1,0)**@{},
(0,4);(1,4)**@{},
(1,0);(2,0)**@{},
(1,4);(2,4)**@{},
(2,0);(3,0)**@{},
(2,4);(3,4)**@{},
(3,0);(4,0)**@{},
(3,4);(4,4)**@{},
(4,0);(5,0)**@{},
(4,4);(5,4)**@{},
(0,0);(0,1)**@{},
(5,0);(5,1)**@{},
(0,1);(0,2)**@{},
(5,1);(5,2)**@{},
(0,2);(0,3)**@{},
(5,2);(5,3)**@{},
(0,3);(0,4)**@{},
(5,3);(5,4)**@{},
\end{xy}
\]
\caption{$nm$ quadrics with $E_4$-points and $2nm$ planes with $E_6$-points}\label{fig:quadrics2}
\end{figure}
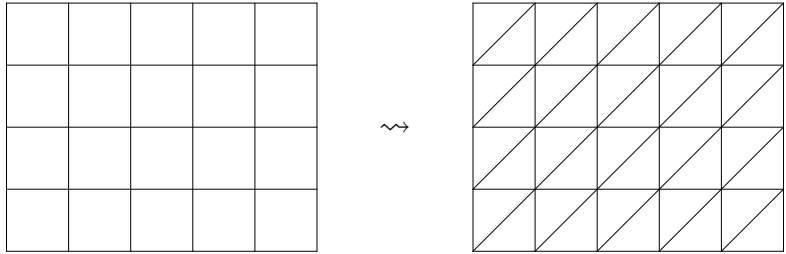

Again each quadric degenerates to the union of two planes. By
doing this as depicted in Figure \ref{fig:quadrics2}, one gets a
degeneration of a general abelian surface with a polarization of
type $(n, m)$ to a planar Zappatic surface of degree $2nm$ with
only $E_6$-points as Zappatic singularities.

Other examples of degenerations, similar to the one considered
above, for $K3$ surfaces (the so called {\em pillow
degenerations}) are considered in e.g. \cite{CMT}.
\end{example}

\begin{remark}
Going back to the general case, if either $C$ or $C'$ has genus greater than 1
and if they degenerate to stick curves, then the surface $S=C\times C'$
degenerates to a union of quadrics, as we said.
Unfortunately it is not clear if it is possible
to further independently degenerate each quadric to two planes.
\end{remark}


From now on, until the end of this section, we deal with degenerations of rational normal scrolls only.
Namely we will show that a general rational normal scroll degenerates to a planar Zappatic
surface with Zappatic singularities of type $R_3$ only and we will see how ``general''
the scroll has to be in order to admit such degenerations
(e.g., in Example \ref{ex:zappa1}, the scrolls are actually forced to have even degree).

There are several ways to construct these degenerations.
We will start from the trivial family and then we will perform
two basic operations: (1) blowing-ups and blowing-downs in the central fibre,
(2) twisting the hyperplane bundle by a component of the central fibre.

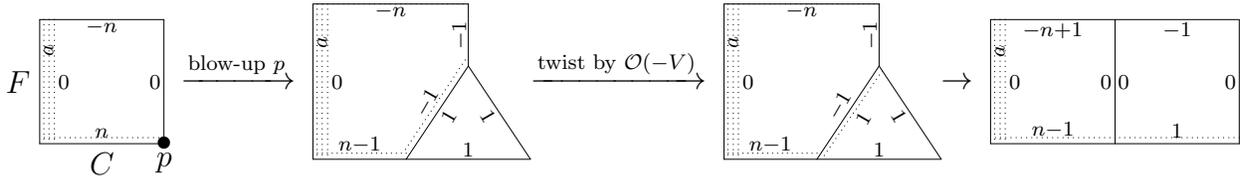
\begin{figure}[ht]
\[
\begin{xy}
0; <4em,0em>:
(0,0);(1,0)*=0{\bullet}**@{-};(1,1)**@{-};(0,1)**@{-};(0,0)**@{-} ,
(0,0);(1,0)**@{} ?*+!U{C} ?(1)*+!U{p},
(0,0)+<0em,0.2em>;(1,0)+<0em,0.2em>**@{.} ?*!D{\st n} ,
(0,0);(0,1)**@{} ?*+!R{F} ,
(0,0)+<0.3em,0em>;(0,1)+<0.3em,0em>**\dir3{.} ?*+!L{\st 0} ?(0.75)*[@!90]{\st a} ,
(0,1);(1,1)**@{} ?*!U{\st -n} ,
(1,0)+<-0.1em,0em>;(1,1)+<-0.1em,0em>**@{} ?*!R{\st 0} ,
\end{xy}
\raisebox{2em}{$\xrightarrow{\st \text{blow-up }p}$}
\;\;
\raisebox{-0.5em}{$%
\begin{xy}
0; <4em,0em>:
(0,0);(0.75,0)**@{-};(1.25,0.75)**@{-};(1.25,1.25)**@{-};(0,1.25)**@{-};(0,0)**@{-} ,
(0.75,0);(1.75,0)**@{-};(1.25,0.75)**@{-} ,
(0,0)+<0em,0.2em>;(0.75,0)+<0em,0.2em>**@{.} ?*!D{\st n-1} ,
(0,0)+<0.3em,0em>;(0,1.25)+<0.3em,0em>**\dir3{.} ?*+!L{\st 0} ?(0.75)*[@!90]{\st a} ,
(0.75,0)+<-0.1em,0.2em>;(1.25,0.75)+<-0.1em,0.2em>**@{.} ?*!D[@!60]{\st -1} ,
(1.25,0.75);(1.25,1.25)**@{} ?*!D[@!90]{\st -1} ,
(0,1.25);(1.25,1.25)**@{.} ?*!U{\st -n} ,
(0.75,0);(1.25,0.75)**@{} ?*!UL[@!60]{\st 1} ,
(0.75,0)+<0em,0.1em>;(1.75,0)+<0em,0.1em>**@{} ?*!D{\st 1},
(1.75,0);(1.25,0.75)**@{} ?*!UR[@!300]{\st 1} ,
\end{xy}%
$}
\raisebox{2em}{$\xrightarrow{\st \text{twist by }\Oc(-V)}$}
\;\;
\raisebox{-0.5em}{$%
\begin{xy}
0; <4em,0em>:
(0,0);(0.75,0)**@{-};(1.25,0.75)**@{-};(1.25,1.25)**@{-};(0,1.25)**@{-};(0,0)**@{-} ,
(0.75,0);(1.75,0)**@{-};(1.25,0.75)**@{-} ,
(0,0)+<0em,0.2em>;(0.75,0)+<0em,0.2em>**@{.} ?*!D{\st n-1} ,
(0,0)+<0.3em,0em>;(0,1.25)+<0.3em,0em>**\dir3{.} ?*+!L{\st 0} ?(0.75)*[@!90]{\st a} ,
(0.75,0)+<-0.05em,0.1em>;(1.25,0.75)+<-0.05em,0.1em>**@{} ?*!D[@!60]{\st -1} ,
(1.25,0.75);(1.25,1.25)**@{} ?*!D[@!90]{\st -1} ,
(0,1.25);(1.25,1.25)**@{.} ?*!U{\st -n} ,
(0.75,0)+<0.15em,0em>;(1.25,0.75)+<0.1em,-0.15em>**@{.} ?*!UL[@!60]{\st 1} ,
(0.75,0)+<0em,0.1em>;(1.75,0)+<0em,0.1em>**@{} ?*!D{\st 1},
(1.75,0);(1.25,0.75)**@{} ?*!UR[@!300]{\st 1} ,
\end{xy}%
$}
\raisebox{2em}{$\to$}
\;\;
\begin{xy}
0; <4em,0em>:
(0,0);(2,0)**@{-};(2,1)**@{-};(0,1)**@{-};(0,0)**@{-} ,
(1,0);(1,1)**@{-} ,
(0,0)+<0em,0.2em>;(1,0)+<0em,0.2em>**@{.} ?*!D{\st n-1} ,
(0,0)+<0.3em,0em>;(0,1)+<0.3em,0em>**\dir3{.} ?*+!L{\st 0} ?(0.75)*[@!90]{\st a} ,
(0,1)+<0em,-0.1em>;(1,1)+<0em,-0.1em> **@{} ?*!U{\st -n+1} ,
(1,0)+<-0.1em,0em>;(1,1)+<-0.1em,0em>**@{} ?*!R{\st 0} ,
(1,0)+<0em,0.2em>;(2,0)+<0em,0.2em>**@{.} ?*!D{\st 1} ,
(1,0)+<0.1em,0em>;(1,1)+<0.1em,0em>**@{} ?*!L{\st 0} ,
(1,1)+<0em,-0.1em>;(2,1)+<0em,-0.1em> **@{} ?*!U{\st -1} ,
(2,0)+<-0.1em,0em>;(2,1)+<-0.1em,0em>**@{} ?*!R{\st 0} ,
\end{xy}
\]
\caption{Degeneration of a scroll $S_{a,b}$ to the union of a plane and a scroll $S_{a,b-1}$} \label{fig:rat1}
\end{figure}

\begin{constr} \label{ex:rat1}
{\em Let $S=S_{a,b}$ be a smooth, rational, normal scroll of type $(a,b)$ in $\pp^r$, where
$r=a+b+1\geq 3$ and we assume that $b\geq a$.
Then $S$ degenerates to the union of a plane and a smooth, rational normal scroll $S_{a,b-1}$
meeting the plane along a ruling.}

\vskip 7pt

\noindent
Indeed, $S$ is the embedding of the Hirzebruch surface $\F_n$, $n=b-a\geq 0$, via the linear system $|C+aF|$,
where $F$ is the ruling and $C$ is a section of self-intersection $n$ (clearly, if $n=0$, we may choose
$F$ to be either one of the two rulings and $C$ to be the other ruling). Set $H=C+aF$.
Consider the trivial family $\s=\F_n\times\D\xrightarrow{\sigma}\D$.
On $\s$ we have the hyperplane bundle $\HH$ which coincides with $H$ on each fibre of $\sigma$.

Now blow up $\s$ at a general point of the central fibre $\s_0$.
Let $V$ be the exceptional divisor and $S'$ be the proper transform of $\s_0$.
Then, $\HH\otimes\Oc(-V)$ embeds $V$ as a plane and maps $S'$ to a scroll of type $(a,b-1)$,
which meet each other along a ruling of $S'$.
We explain these operations in Figure \ref{fig:rat1}, where the dotted lines represent the hyperplane bundle.
The last arrow is the so-called {\em type I transformation} on the vertical $(-1)$-curve (cf.\ \cite{FM}),
which consists in blowing up the $(-1)$-curve and then blowing down the exceptional divisor, which is a $\F_0$, along the other ruling.
The total effect on $\s_0$ is to perform an elementary transformation.

When $r=3$ this process gives the degeneration of a smooth quadric to two planes meeting along a line.
\end{constr}

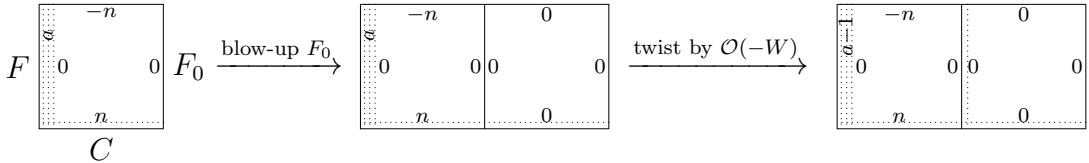
\begin{figure}[ht]
\[
\begin{xy}
0; <4em,0em>:
(0,0);(1,0)**@{-};(1,1)**@{-};(0,1)**@{-};(0,0)**@{-} ,
(0,0);(1,0)**@{} ?*+!U{C} ,
(0,0)+<0em,0.2em>;(1,0)+<0em,0.2em>**@{.} ?*!D{\st n} ,
(0,0);(0,1)**@{} ?*+!R{F} ,
(0,0)+<0.3em,0em>;(0,1)+<0.3em,0em>**\dir3{.} ?*+!L{\st 0} ?(0.75)*[@!90]{\st a} ,
(0,1);(1,1)**@{} ?*!U{\st -n} ,
(1,0)+<-0.1em,0em>;(1,1)+<-0.1em,0em>**@{} ?*!R{\st 0} ,
(1,0);(1,1)**@{} ?*+!L{F_0} ,
\end{xy}
\raisebox{2em}{$\xrightarrow{\st \text{blow-up }F_0}$}
\;\;
\begin{xy}
0; <4em,0em>:
(0,0);(2,0)**@{-};(2,1)**@{-};(0,1)**@{-};(0,0)**@{-} ,
(1,0);(1,1)**@{-} ,
(0,0)+<0em,0.2em>;(1,0)+<0em,0.2em>**@{.} ?*!D{\st n} ,
(0,0)+<0.3em,0em>;(0,1)+<0.3em,0em>**\dir3{.} ?*+!L{\st 0} ?(0.75)*[@!90]{\st a} ,
(0,1)+<0em,-0.1em>;(1,1)+<0em,-0.1em> **@{} ?*!U{\st -n} ,
(1,0)+<-0.1em,0em>;(1,1)+<-0.1em,0em>**@{} ?*!R{\st 0} ,
(1,0)+<0em,0.2em>;(2,0)+<0em,0.2em>**@{.} ?*!D{\st 0} ,
(1,0)+<0.1em,0em>;(1,1)+<0.1em,0em>**@{} ?*!L{\st 0} ,
(1,1)+<0em,-0.1em>;(2,1)+<0em,-0.1em> **@{} ?*!U{\st 0} ,
(2,0)+<-0.1em,0em>;(2,1)+<-0.1em,0em>**@{} ?*!R{\st 0} ,
\end{xy}
\;\;
\raisebox{2em}{$\xrightarrow{\st \text{twist by }\Oc(-W)}$}
\;\;
\begin{xy}
0; <4em,0em>:
(0,0);(2,0)**@{-};(2,1)**@{-};(0,1)**@{-};(0,0)**@{-} ,
(1,0);(1,1)**@{-} ,
(0,0)+<0em,0.2em>;(1,0)+<0em,0.2em>**@{.} ?*!D{\st n} ,
(0,0)+<0.3em,0em>;(0,1)+<0.3em,0em>**\dir3{.} ?*+!L{\st 0} ?(0.75)*[@!90]{\st a-1} ,
(0,1)+<0em,-0.1em>;(1,1)+<0em,-0.1em> **@{} ?*!U{\st -n} ,
(1,0)+<-0.1em,0em>;(1,1)+<-0.1em,0em>**@{} ?*!R{\st 0} ,
(1,0)+<0em,0.2em>;(2,0)+<0em,0.2em>**@{.} ?*!D{\st 0} ,
(1,0)+<0.2em,0em>;(1,1)+<0.2em,0em>**@{.} ?*!L{\st 0} ,
(1,1)+<0em,-0.1em>;(2,1)+<0em,-0.1em> **@{} ?*!U{\st 0} ,
(2,0)+<-0.1em,0em>;(2,1)+<-0.1em,0em>**@{} ?*!R{\st 0} ,
\end{xy}
\]
\caption{Degeneration of a scroll $S_{a,b}$ to the union of a quadric and a scroll $S_{a-1,b-1}$} \label{fig:rat2}
\end{figure}

\begin{constr} \label{ex:rat2}
{\em Let $S=S_{a,b}$ be a smooth, rational, normal scroll of type $(a,b)$ in $\pp^r$, where $r=a+b+1$ and
assume that $b\geq a>1$.
Then $S$ degenerates to the union of a quadric and a smooth, rational normal scroll $S_{a-1,b-1}$
meeting the quadric along a ruling.}

\vskip 7pt

\noindent
Indeed, consider the Hirzebruch surface $\F_n$, $n=b-a\geq 0$, and the trivial family
$\s=\F_n\times\D\xrightarrow{\sigma}\D$,
with the hyperplane bundle $\HH$, as in Construction \ref{ex:rat1}.

Now blow up a ruling $F_0$ in the central fibre $\s_0$.
Let $W$ be the exceptional divisor and $S'$ be the proper transform of $\s_0$.
Then $\HH\otimes\Oc(-W)$ embeds $W$ as a quadric and $S'$ as a scroll of type $(a-1,b-1)$,
which meet along a ruling of $S'$, cf.\ Figure \ref{fig:rat2}.
\end{constr}

By induction on the degree of the scroll and by using Constructions \ref{ex:rat1} and \ref{ex:rat2} for
the inductive steps, we now show the following:

\begin{proposition} \label{prop:rat1}
Let $d\geq 2$ and set $r=d+1\geq3$.
Let $X:= X_{d,0}$ be a planar Zappatic surface of degree $d$ in $\pp^{r}$,
whose dual graph is a chain,
i.e.\ $X$ has $d-2$ points of type $R_3$ as Zappatic singularities.
Then, the Hilbert point of $X$ belongs to the irreducible component ${\mathcal H}_{d,0}$ of
the Hilbert scheme parametrizing rational normal scrolls of degree $d$.
\end{proposition}

\begin{remark} It is well-known (cf.\ e.g. Lemma 3 in \cite{CLM})
that ${\mathcal H}_{d,0}$ is generically reduced and of dimension $d^2+4d-3$.
\end{remark}

\begin{figure}[ht]
\[
\begin{xy}
0; <4em,0em>:
(0,0);(1,0)**@{-};(1,1)**@{-};(0,1)**@{-};(0,0)**@{-} ,
(0,0)+<0em,0.2em>;(1,0)+<0em,0.2em>**@{.} ?*!D{\st 1} , 
(0,0)+<0.3em,0em>;(0,1)+<0.3em,0em>**\dir3{.} ?*+!L{\st 0} ?(0.75)*[@!90]{\st a} ,
(0,1)+<0em,-0.1em>;(1,1)+<0em,-0.1em> **@{} ?*!U{\st -1} ,
(1,0)+<-0.1em,0em>;(1,1)+<-0.1em,0em>**@{} ?*!R{\st 0} ,
\end{xy}
\;\;
\raisebox{2em}{$\to$}
\;\;
\begin{xy}
0; <4em,0em>:
(0,0)*=0{\bullet};(2,0)**@{-};(2,1)**@{-};(0,1)**@{-};(0,0)**@{-} ,
(1,0);(1,1)**@{-} ,
(0,0)+<0em,0.2em>;(1,0)+<0em,0.2em>**@{.} ?*!D{\st 0} , 
(0,0)+<0.3em,0em>;(0,1)+<0.3em,0em>**\dir3{.} ?*+!L{\st 0} ?(0.75)*[@!90]{\st a} ,
(0,1)+<0em,-0.1em>;(1,1)+<0em,-0.1em> **@{} ?*!U{\st 0} ,
(1,0)+<-0.1em,0em>;(1,1)+<-0.1em,0em>**@{} ?*!R{\st 0} ,
(1,0)+<0em,0.2em>;(2,0)+<0em,0.2em>**@{.} ?*!D{\st 1} ,
(1,0)+<0.1em,0em>;(1,1)+<0.1em,0em>**@{} ?*!L{\st 0} ,
(1,1)+<0em,-0.1em>;(2,1)+<0em,-0.1em> **@{} ?*!U{\st -1} ,
(2,0)+<-0.1em,0em>;(2,1)+<-0.1em,0em>**@{} ?*!R{\st 0} ,
\end{xy}
\;\;
\raisebox{2em}{$\to$}
\;\;
\begin{xy}
0; <4em,0em>:
(0,0);(3,0)**@{-};(3,1)**@{-};(0,1)**@{-};(0,0)**@{-} ,
(1,0);(1,1)**@{-} ,
(2,0);(2,1)**@{-} ,
(0,0)+<0em,0.2em>;(1,0)+<0em,0.2em>**@{.} ?*!D{\st 1} ,
(0,0)+<0.3em,0em>;(0,1)+<0.3em,0em>**\dir3{.} ?*+!L{\st 0} ?(0.75)*[@!90]{\st a-1} ,
(0,1)+<0em,-0.1em>;(1,1)+<0em,-0.1em> **@{} ?*!U{\st -1} ,
(1,0)+<-0.1em,0em>;(1,1)+<-0.1em,0em>**@{} ?*!R{\st 0} ,
(1,0)+<0em,0.1em>;(2,0)+<0em,0.1em>**@{} ?*!D{\st -1} ,
(1,0)+<0.1em,0em>;(1,1)+<0.1em,0em>**@{} ?*!L{\st 0} ,
(1,1)+<0em,-0.2em>;(2,1)+<0em,-0.2em> **@{.} ?*!U{\st 1} ,
(2,0)+<-0.1em,0em>;(2,1)+<-0.1em,0em>**@{} ?*!R{\st 0} ,
(2,0)+<0em,0.2em>;(3,0)+<0em,0.2em>**@{.} ?*!D{\st 1} ,
(2,0)+<0.1em,0em>;(2,1)+<0.1em,0em>**@{} ?*!L{\st 0} ,
(2,1)+<0em,-0.1em>;(3,1)+<0em,-0.1em> **@{} ?*!U{\st -1} ,
(3,0)+<-0.1em,0em>;(3,1)+<-0.1em,0em>**@{} ?*!R{\st 0} ,
\end{xy}
\]
\caption{Degeneration of $S_{a,a+1}$ to a planar Zappatic surface with only $R_3$-points} \label{fig:rat3}
\end{figure}
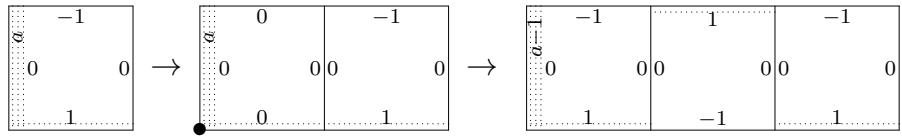

\begin{proof}[Proof of Proposition \ref{prop:rat1}]
We will directly show that a smooth, balanced scroll $S$ degenerates to $X$.

Suppose first that $r$ is even.
Let $S=S(a,a+1)$ be a balanced scroll of degree $d$ in $\pp^r$, i.e.\ $a=(d-1)/2=r/2-1$.
Consider the trivial family $\F_1\times\Delta$, where $\F_1$ is embedded in $\pp^r$ by the linear system $|C+aF|$,
such as in Constructions \ref{ex:rat1} and \ref{ex:rat2}, cf.\ the picture on the left in Figure \ref{fig:rat3}.

Now blow up a ruling in the central fibre, call $W \cong \F_0$ the exceptional divisor
and twist the hyperplane bundle by  $\Oc(-aW)$.
In this way, one gets a degeneration of $S$ to the union of a scroll of type $(a,a)$ in $\pp^{r-1}$
and a plane, meeting along a ruling,
cf.\ Construction \ref{ex:rat2} and the picture in the middle of Figure \ref{fig:rat3}.

Then blow up a general point (the bottom left corner in Figure \ref{fig:rat3}) of the scroll, twist again
by the opposite of the new surface and perform a type I transformation, as we did in Construction \ref{ex:rat1}.
By twisting again by the opposite of the
new surface, counted with multiplicity $a-1$, one gets the configuration depicted
on the right in Figure \ref{fig:rat3},
namely the first two components are two planes, whereas the new component is a scroll of type $(a-1,a)$.

Going on by induction on $a$, by following the same process,
one gets a chain of planes which is a planar Zappatic surface with only $R_3$-points, as wanted.

If $r$ is odd, one starts from a $\F_0$ as in the central picture of Figure \ref{fig:rat3}
and one may perform exactly the same operations in order to get a similar degeneration.
\end{proof}

\begin{remark}\label{rem:fico}
In practice, Proposition \ref{prop:rat1} follows by Contructions
\ref{ex:rat1} and \ref{ex:rat2} with a suitable induction. The
explicit argument we made in the proof shows that there exists a
flat degeneration of smooth, rational scrolls to $X$ whose total
space is singular only at the $R_3$-points of $X$. For another
approach, the reader is also referred to \cite{Mana}.
\end{remark}

\begin{remark}\label{rem:fasci} Suppose to have a smooth scroll $S$ which is the general fibre of
an embedded degeneration in $\Pp^r$ to a Zappatic planar surface $X$. The ruling of $S$, considered
as a curve $\Gamma$ in the Grasmannian $\G(1,r)$, accordingly degenerates to a stick-curve $\Gamma_0$.
This means that
the ruling degenerates to a union of pencils of lines, one in each plane of $X$. Since $\Gamma_0$ is
connected, each double line of $X$ belongs to the pencil in either one of the two planes containing
it.  Hence, the centers of the pencils also belong to the double lines of $X$. Therefore, on each plane
which contains more than one double line of $X$, all the double lines pass through the same
Zappatic singularity which is the center of the pencil. However, the location of the centers of the pencils
on the planes containing only one double line of $X$ is not predictable.
\end{remark}

We conclude this section by proving the following:

\begin{proposition} \label{prop:rat2}
Let $S=S_{a,b}$ be a smooth, rational normal scroll in $\pp^{a+b+1}$, with $b-a\geq4$.
Assume that $S$ is the general fibre of a degeneration whose central fibre is
a planar Zappatic surface $X$.
Then $X$ has worse singularities than $R_3$-points.
\end{proposition}

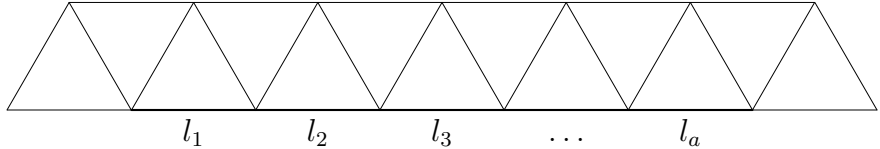
\begin{figure}[ht]
\[
\begin{xy}
0; <4em,0em>:
(0,0);(7,0)**@{-};(6.5,\halfrootthree)**@{-};(0.5,\halfrootthree)**@{-};(0,0)**@{-} ,
(6.5,\halfrootthree);(6,0)**@{-};(5.5,\halfrootthree)**@{-};(5,0)**@{-};
(4.5,\halfrootthree)**@{-};(4,0)**@{-};(3.5,\halfrootthree)**@{-};(3,0)**@{-};
(2.5,\halfrootthree)**@{-};(2,0)**@{-};(1.5,\halfrootthree)**@{-};(1,0)**@{-};
(0.5,\halfrootthree)**@{-} ,
(1,0);(2,0)**[thicker]@{-} ?*+!U{l_1} ,
(2,0);(3,0)**[thicker]@{-} ?*+!U{l_2} ,
(3,0);(4,0)**[thicker]@{-} ?*+!U{l_3} ,
(4,0);(5,0)**[thicker]@{-} ?*+++!U{\ldots} ,
(5,0);(6,0)**[thicker]@{-} ?*+!U{l_a} ,
\end{xy}
\]
\caption{Degeneration of $S_{a,b}$, $b=a+3$, to $X$ with only $R_3$-points} \label{fig:b=a+3}
\end{figure}

\begin{proof}
By construction of the scroll $S$ (cf.\ Definition \ref{def:scroll}), the minimum
degree of a section of $S$ is $a$ and let $C_a$ be the section of degree $a$.
Suppose by contradiction that $S$ is the general fibre of an embedded
degeneration of surfaces whose central fibre is
a planar Zappatic surface $X=\bigcup_{i=1}^{a+b} V_i$ in $\pp^{a+b+1}$, with only
$R_3$-points as Zappatic singularities.
Then the dual graph $G_X$ is a chain and we may and will assume that two planes
$V_i$ and $V_j$ meet along a line if and only if $j=i\pm1$.

While $S$ degenerates to $X$, the ruling of $S$ degenerates to a pencil of lines
$\Lambda_i$ on each plane $V_i$, $i=1,\ldots,a+b$ (cf.\ Remark \ref{rem:fasci})
and the section $C_a$ degenerates
to a chain of lines $l_1,\ldots,l_a$, with $l_i\subset V_{j_i}$, $i=1,\ldots,a$, and we may and will assume
that $j_1<j_2<\cdots<j_a$.

The pencil $\Lambda_1$ has to meet $\bigcup_{i=1}^a l_i$,
hence $V_1$ has to have non-empty intersection with $V_{j_1}$, therefore the
assumption that $X$ has at most $R_3$-points implies that $j_1\leq 3$.
For each $k=2,\ldots,a$, the lines $l_k$ and $l_{k-1}$ meet at a point, so the same argument
implies that $j_k\leq j_{k-1}+2$ (cf.\ Figure \ref{fig:b=a+3}).
It follows that $j_a\leq j_1+2(a-1)\leq 2a+1$.

On the other hand, the pencil $\Lambda_{a+b}$ has to meet $\bigcup_{i=1}^a l_i$, hence $j_a\geq a+b-2$.
In conclusion, one has that:
\[
a+b-2\leq j_a \leq 2a+1,
\]
which contradicts the assumption that $b\geq a+4$.
\end{proof}

For another approach to degenerations of rational scrolls to
unions of planes, the reader is referred to \cite{Mana}.

\begin{remark} By following the lines of the proof of Proposition \ref{prop:rat1} it
is possible to prove that, given $a, b$ positive integers such that $ 0 \leq b-a \leq 3$,  there
exist degenerations whose general fibre is a scroll of type $S(a,b)$ and whose central fibre
is a planar Zappatic surface with only $R_3$-points as Zappatic singularities (cf.\ Figure \ref{fig:b=a+3}).
We will not dwell on this here.
\end{remark}

\section{Degenerations of scrolls: inductive constructions}\label{S:4}

In this section we produce families of smooth scrolls of any genus
$g\geq0$ which degenerate to planar Zappatic surfaces with
Zappatic singularities of types $R_3$ and $S_4$ only.

We start by describing the planar Zappatic surfaces which will be
the limits of our scrolls. We will construct these Zappatic
surfaces by induction on $g$. From now on in this section, we will denote
by $X_{d,g}$ a planar Zappatic surface consisting of $d$ planes and whose
sectional genus is $g$.

We start with the case $g=1$.

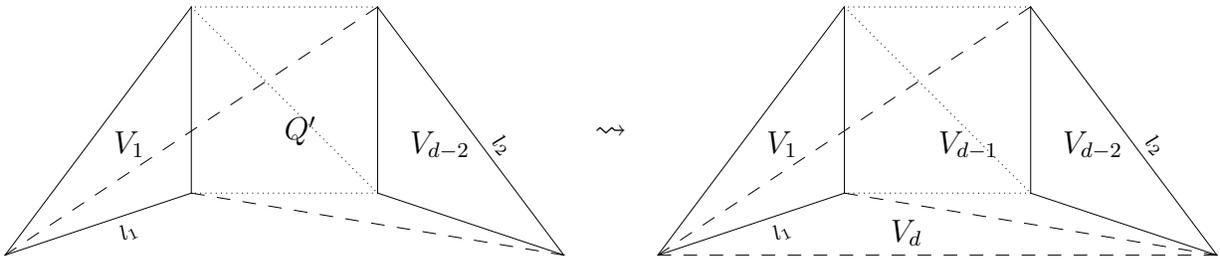
\begin{figure}[ht]
\[
\begin{xy}
0; <4em,0em>: (0,0);(1.5,0.5)**@{-};(1.5,2)**@{-};(0,0)**@{-} ,
(3,0.5);(4.5,0)**@{-};(3,2)**@{-};(3,0.5)**@{-} ,
(1.5,0.5);(3,0.5)**@{.};(1.5,2)**@{.};(3,2)**@{.} ,
(0,0);(1.5,0.5)**@{} ?(.667)*[@!25]+!U{\st l_1} ,
(3,2);(4.5,0)**@{} ?(.667)*[@!315]++!D{\st l_2} ,
(0,0);(3,2)**@{--} , (1.5,0.5);(4.5,0)**@{--} , (1,1)*!U{V_1} ,
(3.5,1)*!U{V_{d-2}} , (2.25,1)*!L{Q'} ,
\end{xy}
\quad \raisebox{4em}{$\rightsquigarrow$} \quad
\begin{xy}
0; <4em,0em>: (0,0);(1.5,0.5)**@{-};(1.5,2)**@{-};(0,0)**@{-} ,
(3,0.5);(4.5,0)**@{-};(3,2)**@{-};(3,0.5)**@{-} ,
(1.5,0.5);(3,0.5)**@{.};(1.5,2)**@{.};(3,2)**@{.} ,
(0,0);(1.5,0.5)**@{} ?(.667)*[@!25]+!U{\st l_1} ,
(3,2);(4.5,0)**@{} ?(.667)*[@!315]++!D{\st l_2} ,
(3,2);(0,0)**@{--};(4.5,0)**@{--};(1.5,0.5)**@{--} , (1,1)*!U{V_1}
, (3.5,1)*!U{V_{d-2}} , (2,0)*+!D{V_d} , (2.5,1)*!U{V_{d-1}} ,
\end{xy}
\]
\caption{Construction of $X_{d,1}$ from $X_{d-2,0}$}
\label{fig:3constr1s}
\end{figure}

\begin{constr}\label{ex:3constr1s}
{\em For any $d\geq 5$, there exists a planar Zappatic surface
$X_{d,1} = \bigcup_{i=1}^d V_i$ in $\Pp^r$, with $r=d-1$, whose
dual graph is a cycle.}

Indeed, if $p_1,\ldots,p_d$ are the coordinate points of $\Pp^r$,
we may let $V_i$, $i=2,\ldots,d-1$, be the plane spanned by $p_{i-1}, p_i, p_{i+1}$
and let $V_1=\langle p_d,p_1,p_2\rangle$, $V_d=\langle p_{d-1},p_d,p_1\rangle$.
Then $X_{d,1} = \bigcup_{i=1}^d V_i$ is a planar Zappatic surface
with dual graph a cycle and whose Zappatic singularities
are points of type $R_3$ at $p_1,\ldots,p_d$, cf.\ Figure \ref{fig:X_d,1},
where one identifies the line $\langle p_d,p_1\rangle$ on the left
with the same line on the right.

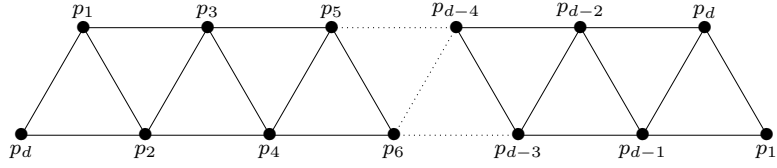
\begin{figure}[ht]
\[
\begin{xy}
0; <4em,0em>:
(0,0)*=0{\bullet};(3,0)**@{-};(4,0)**@{.};(6,0)*=0{\bullet}**@{-};
(5.5,\halfrootthree)*=0{\bullet}**@{-};(3.5,\halfrootthree)**@{-};(2.5,\halfrootthree)**@{.};(0.5,\halfrootthree)*=0{\bullet}**@{-};(0,0)**@{-} ,
(5.5,\halfrootthree);(5,0)*=0{\bullet}**@{-};(4.5,\halfrootthree)*=0{\bullet}**@{-};
(4,0)*=0{\bullet}**@{-};(3.5,\halfrootthree)*=0{\bullet}**@{-};(3,0)*=0{\bullet}**@{.};(2.5,\halfrootthree)*=0{\bullet}**@{-};(2,0)*=0{\bullet}**@{-};
(1.5,\halfrootthree)*=0{\bullet}**@{-};(1,0)*=0{\bullet}**@{-};(0.5,\halfrootthree)**@{-} ,
(0,0)*+!U{\st p_d} ,
(0.5,\halfrootthree)*+!D{\st p_1} ,
(1,0)*+!U{\st p_2} ,
(1.5,\halfrootthree)*+!D{\st p_3} ,
(2,0)*+!U{\st p_4} ,
(2.5,\halfrootthree)*+!D{\st p_5} ,
(3,0)*+!U{\st p_6} ,
(3.5,\halfrootthree)*+!D{\st p_{d-4}} ,
(4,0)*+!U{\st p_{d-3}} ,
(4.5,\halfrootthree)*+!D{\st p_{d-2}} ,
(5,0)*+!U{\st p_{d-1}} ,
(5.5,\halfrootthree)*+!D{\st p_{d}} ,
(6,0)*+!U{\st p_1} ,
\end{xy}
\]
\caption{Planar Zappatic surface $X_{d,1}$ with dual graph a cycle} \label{fig:X_d,1}
\end{figure}

We will show in Theorem \ref{thm:deg3} that $X_{d,1}$ is the flat
limit of a smooth scroll of genus 1 in $\Pp^r$. In order to do
that, now we describe another way to construct $X_{d,1}$, which
will also help to understand the next inductive steps.

Let $X_{d-2,0}=\bigcup_{i=1}^{d-2} V_i$ be a planar Zappatic surface of degree $d-2$ in
$\Pp^r$, whose dual graph is a chain. We may and will assume that
the planes $V_i$ and $V_j$ meet along a line
if and only if $j = i \pm 1$.

Now choose a general line $l_1$ in $V_1$ and a general line $l_2$
in $V_{d-2}$, thus $l_1$ [resp.\ $l_2$] does not pass through the
$R_3$-point $V_1\cap V_2\cap V_3$ [resp.\ $V_{d-4}\cap V_{d-3}\cap
V_{d-2}$]. Clearly the lines $l_1$ and $l_2$ are skew and span a
$\pp^3$, call it $\Pi$. By a computation in coordinates one proves
that, if $d \geq 6$, then $\Pi \cap X_0 = l_1 \cup l_2$. Therefore there exists a
smooth quadric $Q'$ in $\Pi$ such that $l_1$, $l_2$ are lines of
the same ruling on $Q'$ and $Q'$ meets  $X_0$ transversally along
$Q'\cap X_0=l_1\cup l_2$. On the other hand, if $d =5$,
then $\Pi \cap X_0 = l_1 \cup l_2 \cup l$, where $l$ is a line in the central plane.
Nonetheless it is still true that there exists a smooth quadric $Q'$ which contains
$l_1$ and $l_2$ and meets $X_0$ transversally.

Finally, in $\Pi$, the quadric $Q'$
degenerates to two planes $V_{d-1}$ and $V_{d}$, such that $l_i
\subset V_{d-i+1}$, $i=1,2$. By construction, the planar Zappatic
surface $X_{d,1} =X_{d-2,0} \cup V_{d-1}\cup
V_{d}=\bigcup_{i=1}^{d} V_i$ has dual graph which is a cycle,
hence it has only $R_3$-points as Zappatic singularities (cf.\
Example \ref{ex:zappa2} and Figure \ref{fig:3constr1s}). Note
that, if $d \geq 6$, then there are pairs of disjoint planes in
the cycle.
\end{constr}

Next, we complete the construction proceeding inductively.

\begin{figure}[ht]
\[
\begin{xy}
0; <4em,0em>:
(0,0);(1.5,0.5)*=0{\bullet}**@{-};(1.5,2)**@{-};(0,0)**@{-} ,
(3,0.5);(4.5,0)**@{-};(3,2)*=0{\bullet}**@{-};(3,0.5)**@{-} ,
(1.5,0.5);(3,0.5)**@{.};(1.5,2)**@{.};(3,2)**@{.} ,
(1.5,0.5)*+!U{p_1} , (3,2)*+!D{p_2} , (1.5,0.5);(0.75,1)**@{--}
?(.67)*[@!330]+!U{\st l_1} , (3,2);(3.75,0.25)**@{--}
?(.75)*[@!300]++!D{\st l_2} , (0.75,1);(3,2)**@{--} ,
(3.75,0.25);(1.5,0.5)**@{--} , (0.5,1)*!C{V_1} , (4,1)*!C{V_2} ,
(2.25,1)*!C{Q'} , (1.5,0.5);(1.5,-0.5)**@{.};(0,0)**@{.} ,
(3,2);(4.5,2)**@{.};(4.5,0)**@{.} ,
\end{xy}
\quad \raisebox{4em}{$\rightsquigarrow$} \quad
\begin{xy}
0; <4em,0em>:
(0,0);(1.5,0.5)*=0{\bullet}**@{-};(1.5,2)**@{-};(0,0)**@{-} ,
(3,0.5);(4.5,0)**@{-};(3,2)*=0{\bullet}**@{-};(3,0.5)**@{-} ,
(1.5,0.5);(3,0.5)**@{.};(1.5,2)**@{.};(3,2)**@{.} ,
(1.5,0.5)*+!U{p_1} , (3,2)*+!D{p_2} , (1.5,0.5);(0.75,1)**@{-}
?(.67)*[@!330]+!U{\st l_1} , (3,2);(3.75,0.25)**@{-}
?(.75)*[@!300]++!D{\st l_2} ,
(3,2);(0.75,1)**@{-};(3.75,0.25)**@{-};(1.5,0.5)**@{-}
?(.67)*+!U{V_d} , (0.5,1)*!C{V_1} , (4,1)*!C{V_2} ,
(2.5,1)*!D{V_{d-1}} , (1.5,0.5);(1.5,-0.5)**@{.};(0,0)**@{.} ,
(3,2);(4.5,2)**@{.};(4.5,0)**@{.} ,
\end{xy}
\]
\caption{Construction of $X_{d,g}$ from $X_{d-2,g-1}$}
\label{fig:3constr}
\end{figure}
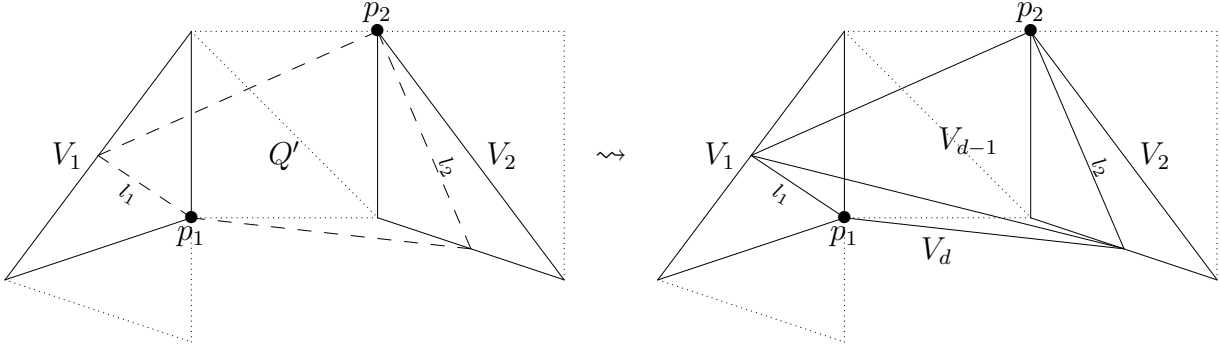

\begin{constr} \label{ex:3constr}
{\em Fix integers $g,d$ such that $g\geq2$ and $d\geq 2g+4$. Set
$c=d-2g-4\geq0$ and $r=5+c=d-2g+1$. There is a planar Zappatic
surface $X_{d,g}=\bigcup_{i=1}^{d} V_i$ in $\Pp^r$ such that:
\begin{itemize}
\item $X_{d,g}$ has $3g+6+c$ double lines, i.e.\ its dual graph
$G_{X_{d,g}}$ has $3g+6+c$ edges; \item $X_{d,g}$ has $r+1$ points
of type $R_3$ and $2(g-1)$ points of type $S_4$; \item for each
$i$, $V_i$ is the central plane through a point $p$ of type either
$R_3$ or $S_4$, i.e.\ $V_i$ is the central component of $X_{d,g}$
passing through $p$ as defined in Notation \ref{nota:central};
\item there exist two $R_3$-points of $X_{d,g}$ whose central
planes do not meet; \item $\chi(\Oc_{X_{d,g}})=1-g$,
$p_\omega(X_{d,g})=0$, $q(X_{d,g})=g(X_{d,g})=g$.
\end{itemize}
}

Taking into account Construction \ref{ex:3constr1s}, which covers
$g=1$ and $d \geq 5$, we can proceed by induction and assume that
we have the surface $X_{d-2, g-1}$.  Let $V_{1}$ and $V_{2}$ be
disjoint planes in $X_{d-2, g-1}$ such that each one of them is
the central plane for a $R_3$-point, say $p_{1}$ and $p_{2}$
respectively.

Now choose a line $l_1$ in $V_{1}$ [resp.\ $l_2$ in $V_{2}$] which
is general among those passing through $p_{1}$ [resp.\ through
$p_{2}$]. Then $l_1$ and $l_2$ are skew and span a $\pp^3$, say
$\Pi$, therefore there exists a smooth quadric $Q'$ in $\Pi$
containing $l_1$ and $l_2$ as lines of the same ruling, cf.\
Figure \ref{fig:3constr}.

Now we prove the following:

\begin{claim}\label{claim:bast}
For general choices, $Q'$ and $X_{d-2, g-1}$ meet transversally
along $X_{d-2, g-1}\cap Q'=l_1\cup l_2$.
\end{claim}
\begin{proof} In order to prove the claim, it suffices to show that $\Pi$ does not meet the remaining components of
$X_{g-1}$ along a curve, i.e.\ that $\Pi$ does not meet $V_i$,
$i\ne 1, 2$, along a line. Before proving the claim, we make a
remark. Suppose that there are two further planes, say $V_{3}$ and
$V_{4}$, in $X_{d-2, g-1}$ contained in $\langle V_{1}, V_{2}
\rangle = \Sigma \cong \pp^5$. Suppose also that  the dual graph
of the planar Zappatic surface $V_{1}\cup V_{3}\cup V_{4}\cup
V_{2}$ is a chain of length 4. Then the points $V_{1}\cap
V_{3}\cap V_{4}$ and $V_{3}\cap V_{4}\cap V_{2}$ are of type
$R_3$. Note that this certainly happens if $c=0$ and $g= 2$
because in that case the dual graph of $X_{d-2, g-1}$ is a cycle
of length six.

In this situation, a computation in coordinates in $\Sigma$ shows
that for a general choice of $l_1$ and $l_2$, $\Pi = \langle l_1,
l_2 \rangle$ does not intersect either $V_{3}$ or $V_{4}$ along a
line.

Now we prove the claim arguing by contradiction. Fix the line
$l_2$ in $V_{2}$ and consider $\langle l_2 , V_{1} \rangle =
\Omega \cong \Pp^4$. By moving $l_1$ in the pencil of lines of
$V_{1}$ through $p_{1}$, one gets a pencil $\Phi$ of $\pp^3$'s
inside $\Omega$ and each of these $\Pp^3$'s meets a plane, say
$V_3$, along a line. There are two possibilities: either $V_3
\subset\Omega$, or $V_3\nsubseteq\Omega$.

In the former case, $V_3$ intersects $V_1$ at a point $q$. Let
$l_2$ move in the pencil of lines of $V_{2}$ through $p_{2}$: one
gets a pencil of $\pp^4$'s in $\Sigma = \langle V_{1}, V_{2}
\rangle$, whose base-locus is $\langle p_1, V_2 \rangle \cong
\Pp^3$ in which $V_3$ is contained. This implies that $q=p_1$,
moreover $V_3$ intersects $V_2$ along a line which necessarily
contains $p_2$. In conclusion, $V_3$ contains the line passing
through $p_1$ and $p_2$. This yields the existence of a plane
$V_4$ which forms, together with $V_1$, $V_2$ and $V_3$, a
configuration in $\Sigma$ of four planes as the one discussed
above. This is a contradiction.

Suppose now that $V_3\nsubseteq\Omega$. Then $V_3$ meets along a
line the base locus of the pencil $\Phi$, which is the plane
$\langle p_1, l_2 \rangle$. By moving $l_2$, we see that $V_3$ has
to contain the line through $p_1$ and $p_2$ and we get a
contradiction as before.
\end{proof}

In $\Pi$, the smooth quadric $Q'$ degenerates to the union of two
planes, say $V_{d-1}\cup V_{d}$, where $l_i\subset V_{d-i+1}$,
$i=1,2$. Consider the planar Zappatic surface $X_{d,g}=X_{d-2,
g-1}\cup V_{d-1}\cup V_{d}$ of degree $d$ in $\pp^{r}$. Thus, we
added to $X_{d-2, g-1}$ two planes and three double lines
$V_{1}\cap V_{d}$, $V_{d}\cap V_{d-1}$ and $V_{d-1}\cap V_{2}$.
Moreover, the points $p_{1}$ and $p_{2}$ become points of type
$S_4$ for $X_{g}$ and we added two further points of type $R_3$ at
$V_{1}\cap V_{d}\cap V_{d-1}$ and $V_{d}\cap V_{d-1} \cap V_{2}$,
cf.\ Figure \ref{fig:3constr}. Finally, one checks that each one
of the planes $V_{d-1}$ and $V_d$ is disjoint from some other
plane in the configuration. This ends the construction.
\end{constr}

Next, we will prove that the Zappatic surfaces $X_{d,g}$ we
constructed are limits of smooth scrolls of genus $g$. First we
make a remark.

\begin{remark}\label{rem:xg}
If $X_{d,g}$ is the flat limit of a family of smooth surfaces $Y$,
then Theorem \ref{thm:invXt} implies that:
\begin{equation} \label{eq:invS}
g(Y) = g,
\qquad
p_g(Y) = 0,
\qquad
\chi(\Oc_{Y})=1-g,
\qquad
8(1-g)\leq K^2_Y \leq 6(1-g).
\end{equation}
\end{remark}

\begin{theorem}\label{thm:deg3} Let $g \geq 0$ and $d \geq 2g+4$ be integers. Let $r = d-2g+1$.
The Hilbert point
corresponding to the planar Zappatic surface $X_{d,g}$ belongs to an irreducible component
${\mathcal H}_{d,g}$ of the Hilbert scheme of scrolls of degree $d$ and genus $g$ in $\Pp^{r}$,
such that:
\begin{itemize}
\item[(i)] the general point of $\HH_{d,g}$ represents a smooth, linearly normal scroll $Y \subset \Pp^r$;
\item[(ii)]  $\HH_{d,g}$ is generically reduced, $dim (\HH_{d,g}) = h^0(Y, \N_{Y/\Pp^r}) = (r+1)^2+7(g-1) $,
and moreover $h^1(Y, \N_{Y/\Pp^r})= h^2(Y, \N_{Y/\Pp^r}) = 0$.
\end{itemize}
\end{theorem}
\begin{proof}[Proof of Theorem \ref{thm:deg3}: beginning]
We prove Theorem \ref{thm:deg3} by induction on $g$. The case $g= 0$ has been treated in Proposition
\ref{prop:rat1}. By induction on $g$, we may assume that
$X_{d-2, g-1}$ is the flat limit of a  smooth scroll $S$
of degree $d-2$ and genus $g-1$ in $\pp^{r}$,
which is represented by a smooth point of a component ${\mathcal H}_{d-2,g-1}$ of
the Hilbert scheme of dimension $(r+1)^2+7(g-2)$.

We can now choose $l_1$ and $l_2$ as in Constructions \ref{ex:3constr1s} and \ref{ex:3constr} so that they
are limits of
rulings $F_1$ and $F_2$, respectively, on $S$ (cf.\ Remark \ref{rem:fasci}).

Let $Q$ be a smooth quadric containing $F_1$ and $F_2$, whose limit is $Q'$.
By the properties of $X_{d-2, g-1}$ and of $Q'$ (see Claim \ref{claim:bast}), it follows that $S$ and $Q$
meet transversally along $S\cap Q=F_1\cup F_2$.

The inductive step is a consequence of the following lemma.
\end{proof}

\begin{lemma}\label{lem:deg3} In the above setting, consider the union
\[
R:= S \cup Q.
\]
Let $\N_R$ and $\T_R$ be the normal and the tangent sheaf of $R$ in $\Pp^r$, respectively; then,  one has:
\begin{gather}
H^1(\N_R) = H^2(\N_R) = 0,            \label{eq:deg3a} \\
h^0(\N_R) = (r+1)^2+7(g-1)=d^2-4dg+4d+4g^2-g-3.  \label{eq:deg3b}
\end{gather}
Furthermore the natural map $H^0(\N_R)\to H^0(T^1)$, induced by the exact sequence
\begin{equation}\label{eq:t1a}
0 \to \T_R \to \T_{\Pp^r}|_R \xrightarrow{\,\tau\,} \N_R \to T^1 := Coker(\tau) \to 0,
\end{equation}is surjective.
\end{lemma}

\begin{proof}
We will compute the cohomology of $\N_R$, by using
a similar technique as in section 2.2 of \cite{CLM} (see Lemma 3 therein).

Let $\Gamma:= S \cap Q=F_1\cup F_2$ be the double curve of $R$. Since $R$ has global normal crossings,
the sheaf $T^1$ in \eqref{eq:t1a} is locally free, of rank 1 on the singular locus $\Gamma$ of $R$
and, by \cite{F}, it is
\[
T^1\cong \N_{\Gamma/S}\otimes \N_{\Gamma/Q}.
\]
Since $\Gamma$ is the union of two lines of the same ruling on both $Q$ and $S$,
it follows that
\begin{equation}\label{eq:t1a2}
T^1 \cong \Oc_{\Gamma}.
\end{equation}
Let us consider the inclusions $\iota_S: \N_S\to\N_R|_S$ and $\iota_Q: \N_Q\to\N_R|_Q$.
Lemma 2 in \cite{CLM} shows that $T^1\cong\coker(\iota_S)$ and $T^1\cong\coker(\iota_Q)$.
For readers' convenience, we recall here the proof.
By a local computation, one sees that the cokernel $K$ of $\iota_S$
is locally free of rank 1 on $\Gamma$.
In the diagram
\begin{equation} \label{eq:diagram}
\xymatrix@R-1.8em{%
\T_{\Pp^r}|_R \ar[r] \ar[dd] & \N_R \ar[r] \ar[rd] & T^1  \ar[r] & 0 \\
                  &       & \N_R|_S  \ar[rd]                      \\
\T_{\Pp^r}|_S \ar[r] & \N_S \ar[ru]^{\iota_S} &              & 0
}
\end{equation}the horizontal and diagonal rows are exact,
hence the commutativity of the pentagon shows that $T^1$ surjects onto $K$.
Since both are locally free sheaves of rank 1, one concludes that $T^1\cong K$.
The same argument works for $Q$.

Hence
the following sequences are exact:
\begin{gather}
0 \to \N_S \to \N_R|_S \to T^1 \to 0,                             \label{eq:T1coker} \\
0 \to \N_Q(-\Gamma) \to \N_R|_Q(-\Gamma) \to T^1(-\Gamma) \to 0.  \label{eq:T1coker2}
\end{gather}

Moreover, one has the exact sequence
\begin{equation}\label{eq:deg3.2}
0 \to \N_R|_Q \otimes \Oc_R(-\Gamma) \to \N_R \to \N_R|_S \to 0,
\end{equation}
so that, in order to prove \eqref{eq:deg3a}, it suffices to show that
\begin{align}
& H^i(\N_R|_S) =  0,                         &\text{for }1 \leq i \leq 2,  \label{eq:deg3a1} \\
& H^i(\N_R|_Q \otimes \Oc_R(-\Gamma)) = 0,   &\text{for }1 \leq i \leq 2.  \label{eq:deg3a2}
\end{align}
By induction on $g$, one knows that $H^i(\N_S)=0$, $i=1,2$.
By \eqref{eq:t1a2}, one has that $H^i(T^1)=H^i(\Oc_\Gamma)=0$, $i=1,2$,
because $\Gamma$ is the union of two distinct lines.
Hence the  sequence \eqref{eq:T1coker} implies \eqref{eq:deg3a1}.

Note that $H^i(T^1(-\Gamma))=H^i(\Oc_{\Gamma}(-\Gamma))=0$, $i=1,2$.
Taking into account the exact sequence \eqref{eq:T1coker2},
the proof of \eqref{eq:deg3a2}  is concluded if one shows that
\begin{equation} \label{eq:deg3.a3}
H^i(\N_Q(-\Gamma)) = 0, \qquad\text{for }1 \leq i \leq 2.
\end{equation}
Since $Q$ lies in a $\pp^3$, one has that
\[
\N_{Q} \cong \Oc_{Q}(2) \oplus \Oc_{Q}(1)^{\oplus (r-3)}.
\]
Recall that $F_1$ and $F_2$ are lines of the same ruling, so $F_1\sim F_2$ and $\Oc_Q(-\Gamma)\cong\Oc_Q(-2F_1)$.
Let $G$ be the other ruling of $Q$ and $H$ be the general hyperplane section of $Q$,
hence $H\sim G+F_1$ and one has that:
\begin{equation}\label{eq:nq}
\N_{Q} \otimes \Oc_{Q} ( - \Gamma) \cong \Oc_{Q}(2 G) \oplus \Oc_{Q}(G - F_1)^{\oplus (r-3)}
\end{equation}
and one sees that $h^i(\Oc_{Q}(2 G))=h^i(\Oc_{Q}(G-F))=0$, for $i=1,\  2$, which proves
\eqref{eq:deg3.a3}. The proof of \eqref{eq:deg3a} is thus concluded.

\bigskip
We now prove formula \eqref{eq:deg3b}.
By \eqref{eq:deg3a}, one has that $h^0(\N_R)=\chi(\N_R)$,
which one  computes by using \eqref{eq:T1coker},  \eqref{eq:T1coker2} and \eqref{eq:deg3.2}:
\[
\chi(\N_R) = \chi(\N_R|_{S})+\chi(\N_R|_{Q}\otimes\Oc_Q(-\Gamma)) = \chi(\N_S)+
\chi(T^1)+\chi(\N_Q(-\Gamma))+\chi(T^1(-\Gamma)).
\]
By \eqref{eq:t1a2}, one has that $\chi(T^1)=\chi(T^1(-\Gamma))=2$.
By \eqref{eq:nq}, one has $\chi(\N_Q(-\Gamma))=3$.
Finally, by induction
\[
\chi(\N_S)=(r+1)^2+7(g-2),
\]
which concludes the proof of \eqref{eq:deg3b}.

It remains to show that the map $H^0(\N_R)\to H^0(T^1)$ is surjective.
Since $H^1(\N_S)= 0$,
the map $H^0(\N_R|_S) \to H^0(T^1)$
is surjective  by \eqref{eq:T1coker}.
Finally \eqref{eq:deg3a2} implies that $H^0(\N_R)$ surjects onto $H^0(\N_R|_S)$,
which concludes the proof of the lemma.
\end{proof}

We are finally ready for the

\begin{proof}[Proof of Theorem \ref{thm:deg3}: conclusion]
By Lemma \ref{lem:deg3}, one has that $H^1(\N_R) = 0$, which means that $R$ corresponds to a smooth point $[R]$
of the Hilbert scheme of surfaces with degree $d$ and sectional genus $g$ in $\Pp^{d-2g+1}$.
Therefore, $[R]$ belongs to a single reduced component
${\mathcal H}_{d,g}$ of the Hilbert scheme of dimension $h^0(\N_R)$.
The last assertion of Lemma \ref{lem:deg3} implies that a general tangent vector to  ${\mathcal H}_{d,g}$ at the point
$[R]$ represents a first-order embedded deformation of $R$ which smooths the double curve $\Gamma$. Therefore, the general
point in ${\mathcal H}_{d,g}$ represents a smooth, irreducible surface $Y$.
Thus $Y$ degenerates to $R$ and also to the planar Zappatic surface $X_{d,g}$
(cf.\ Proposition \ref{prop:rat1} and Constructions \ref{ex:3constr1s}, \ref{ex:3constr}).

Classical adjunction theory (cf.\ e.g.
\cite{Ionescu} and \S\ 7 in \cite{CR}) implies that $Y$ is a scroll: otherwise, if $H$ is the hyperplane section of $Y$,
one has $K_Y + H$ nef and therefore $ 0 < d \leq 4(g-1) + K^2_Y $ contradicting $K_Y^2 \leq 6(1-g)$ in
\eqref{eq:invS}.

Finally, the assertion about linear normality is trivial for $g= 0$ and is clear by induction and construction,
for $g >0$.
\end{proof}

\begin{remark}\label{rem:ex6}
By using the same first part of the proof of Theorem \ref{thm:deg3}, one can
observe that Construction \ref{ex:3constr} can be carried on also when $d=2g+3$.

Indeed, in this case, $X_{d,g}$ is a union of planes lying in $\pp^4$ which
is not a Zappatic surface if $g\geq2$, since there are singular points
where only two planes of the configuration meet, which are not Zappatic singularities.
The only difference in the construction is that,
since there are no pairs of disjoint planes,
we have to choose $l_1$ and $l_2$ on two planes $V_1$ and $V_2$ which
meet at a point but not along a line.
Moreover the proof of the existence of the quadric meeting transversally
the union of planes along $l_1\cup l_2$ is a bit more involved.

Nonetheless, as in the proof of Theorem
\ref{thm:deg3}, one can show that $X_{d,g}$
is a flat limit of a family of linearly normal scrolls in $\pp^4$ for
any genus $g\geq0$ and degree $d=2g+3$.
These scrolls are smooth only if $g=0,1$,
whereas they have isolated double points if $g\geq2$.
\end{remark}

We finish this section by mentioning two more examples of configurations of planes forming a planar Zappatic surface,
with only points of type $R_3$ and $S_4$, which are degenerations of smooth scrolls.
The advantage of this construction is that they are slightly simpler than Construction \ref{ex:3constr}.
The disadvantage is that they work only for larger values of the degree.

\begin{example}
{\em Fix arbitrary integers $g,d$ such that $g\geq2$ and $d > 4g$. Set
$r=d-2g+1$.
Let $X_{d-2g,0}=\bigcup_{i=1}^{d-2g} V_i$ be a planar Zappatic surface
in $\pp^r$ whose dual graph is a chain.
One can attach $2g$ planes to $X_{d-2g,0}$ in order to get a planar
Zappatic surface $Y_{d,g}$
of degree $d$ and sectional genus $g$ in $\pp^r$ with $d-2g+2$ points of
type $R_3$ and $2g-2$ points of type $S_4$.
}

\vskip 6pt
Indeed, we may assume that $V_i$ meets $V_j$ along a line if and
only if $j=i\pm1$.
Denote by $p_2,\ldots,p_{d-2g-1}$ the points of type $R_3$ of
$X_{d-2g,0}$,
where $p_i=V_{i-1}\cap V_{i}\cap V_{i+1}$, $i=2,\ldots,d-2g-1$.

Choose a general line $l_{1,1}$ in $V_1$ [resp.\ $l_{1,2}$ in $V_{d-2g}$], i.e. a line
not passing through $p_2$ [resp.\ $p_{d-2g-1}$].
For $i=2,\ldots,g$, choose a line $l_{i,1}$ in $V_i$ [resp.\ a line
$l_{i,2}$ in $V_{d-2g+1 -i }$],
which is general among those lines passing through $p_i$ [resp.\ through
$p_{d-2g+1 -i}$].

The generality assumption implies that all the lines $l_{i,1},l_{i,2}$, $ 1 \leq i \leq g$,
are pairwise skew. For every $i=1,\ldots,g$, there is a smooth quadric surface $Q'_i$ which
contains $l_{i,1}$ and $l_{i,2}$, in the $\pp^3$ spanned by them.
In this $\pp^3$ the quadric $Q'_i$ degenerates to two distinct planes,
say $V_{i,1}$ and $V_{i,2}$, leaving $l_{i,1}$ and $l_{i,2}$ fixed: the plane
$V_{i,1}$ contains $l_{i,1}$ whereas $V_{i,2}$ contains $l_{i,2}$.
Then $Y= Y_{d,g}:= X_{d-2g,0 } \cup \bigcup_{i=1}^{g} (V_{i,1} \cup V_{i,2})$ is a planar Zappatic surface in
$\pp^r$. Note that we added to the points $p_2, \ldots , p_{d-2g -1}$ new Zappatic singularities
at the points:
\begin{itemize}
\item[(i)] $q_{i,j}$, with $  1\leq i \leq g$, $1 \leq j \leq 2$, where
$q_{i,1} = V_i \cap V_{i,1} \cap V_{i,2}$ and $q_{i,2} = V_{i,1} \cap V_{i,2} \cap V_{d-2g + 1 -i}$,
\item[(ii)] $p_1 = V_{1} \cap V_{2} \cap V_{1,1}$ and
$p_{d-2g} = V_{d-2g} \cap V_{d-2g-1} \cap V_{1,2}$
\end{itemize}Then $Y$ is a planar Zappatic surface with the following properties:
\begin{itemize}
\item the dual graph $G_Y$ has $d$ vertices and $d+g-1$ edges;
\item $Y$ has $2g-2$ points of type $S_4$, namely
$p_2, \ldots,p_{g}, p_{d-3g + 1}, \ldots , p_{d-2g -1}$;
\item $Y$ has $d-2g+2$ points of type $R_3$, namely
$q_{i,j}$, $1 \leq i \leq g$, $ 1 \leq j \leq 2$, $p_1$, $p_{d-2g}$
and $p_{g+1}, \ldots, p_{d-3g}$.
\item $\chi(\Oc_X)=1-g$, $p_\omega(X)=0$, $q(X)=g(X)=g$,
\end{itemize}(cf. Figure \ref{fig:d>4g}).

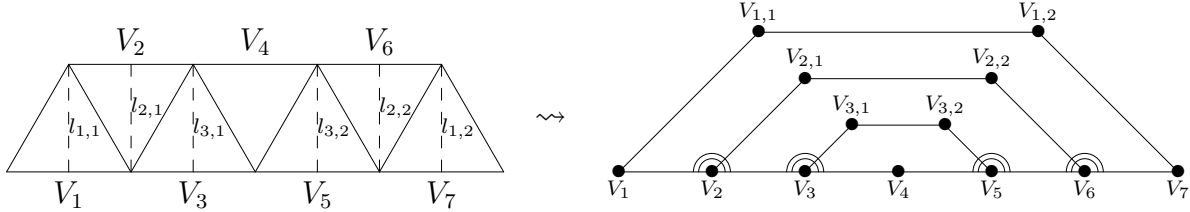
\begin{figure}[ht]
\[
\begin{xy}
0; <4em,0em>:
(0,0);(4,0)**@{-};(3.5,\halfrootthree)**@{-};(0.5,\halfrootthree)**@{-};(0,0)**@{-} ,
(3.5,\halfrootthree);(3,0)**@{-};(2.5,\halfrootthree)**@{-};(2,0)**@{-};
(1.5,\halfrootthree)**@{-};(1,0)**@{-};(0.5,\halfrootthree)**@{-} ,
(0.5,\halfrootthree);(0.5,0)**@{--} ?*!UL{\st l_{1,1}} ?(1)*+!U{V_1},
(1,0);(1,\halfrootthree)**@{--} ?*!DL{\st l_{2,1}} ?(1)*+!D{V_2},
(1.5,\halfrootthree);(1.5,0)**@{--} ?*!UL{\st l_{3,1}} ?(1)*+!U{V_3},
(2,0);(2,\halfrootthree)**@{} ?(1)*+!D{V_4},
(2.5,\halfrootthree);(2.5,0)**@{--} ?*!UL{\st l_{3,2}} ?(1)*+!U{V_5},
(3,0);(3,\halfrootthree)**@{--} ?*!DL{\st l_{2,2}} ?(1)*+!D{V_6},
(3.5,\halfrootthree);(3.5,0)**@{--} ?*!UL{\st l_{1,2}} ?(1)*+!U{V_7},
\end{xy}
\quad
\raisebox{1.75em}{$\rightsquigarrow$}
\quad
\begin{xy}
0; <3em,0em>:
(0,0);(6,0)**@{-} ,
(0,0)*=0{\bullet};(1.5,1.5)*=0{\bullet}**@{-};(4.5,1.5)*=0{\bullet}**@{-};(6,0)*=0{\bullet}**@{-} ,
(1,0)*=0{\bullet};(2,1)*=0{\bullet}**@{-};(4,1)*=0{\bullet}**@{-};(5,0)*=0{\bullet}**@{-} ,
(2,0)*=0{\bullet};(2.5,0.5)*=0{\bullet}**@{-};(3.5,0.5)*=0{\bullet}**@{-};(4,0)*=0{\bullet}**@{-} ,
(3,0)*=0{\bullet} ,
(0,0)*+!U{\st V_1} ,
(1,0)*+!U{\st V_2}*\cir<5pt>{u^d}*\cir<7pt>{u^d} ,
(2,0)*+!U{\st V_3}*\cir<5pt>{u^d}*\cir<7pt>{u^d} ,
(3,0)*+!U{\st V_4} ,
(4,0)*+!U{\st V_5}*\cir<5pt>{u^d}*\cir<7pt>{u^d} ,
(5,0)*+!U{\st V_6}*\cir<5pt>{u^d}*\cir<7pt>{u^d} ,
(6,0)*+!U{\st V_7} ,
(1.5,1.5)*+!D{\st V_{1,1}} ,
(4.5,1.5)*+!D{\st V_{1,2}} ,
(2,1)*+!D{\st V_{2,1}} ,
(4,1)*+!D{\st V_{2,2}} ,
(2.5,0.5)*+!D{\st V_{3,1}} ,
(3.5,0.5)*+!D{\st V_{3,2}} ,
\end{xy}
\]
\caption{Construction of $Y_{d,g}$ from $X_{d-2g,0}$ for $g=3$ and $d=4g+1=13$} \label{fig:d>4g}
\end{figure}

Recall that $X_{d-2g,0}$ is the flat limit of a smooth, rational normal
scroll $S$ of degree $d-2g$ in $\Pp^{d-2g+1}$. If
$F_{i,j}$, $ 1 \leq i \leq g$, $1 \leq j \leq 2$, is
the ruling of $S$ whose limit is $l_{i,j}$ and $Q_i$
a smooth quadric containing $F_{i,1}$, $F_{i,2}$, whose limit is $Q'_i$,
then one can show, by using similar techniques as in the proof of
Theorem \ref{thm:deg3},
that the union of the rational normal scroll $S$ and the $g$ quadrics
$Q_i$ is
a flat limit of a family of smooth, linearly normal scrolls of degree $d$ and
genus $g$ in $\Pp^{d-2g+1}$, which is contained in a the same component
${\mathcal H}_{d,g}$ of Theorem \ref{thm:deg3} (cf.\ Theorem \ref{thm:contoflam} and
Remark \ref{rem:new} below).
\end{example}

With a slight modification of the previous construction, one can cover also the case $d=4g$.
We do not dwell on this here.

\begin{example} {\em Fix integers $g,d$ such that $g\geq1$ and $d\geq 3g+2$.
By induction on $g$, we will construct a planar Zappatic surface
$Z_{d,g}=\bigcup_{i=1}^{d} V_i$
in $\pp^{d-2g+1}$ such that:
\begin{itemize}
\item $Z_{d,g}$ has $d-2g+1$ double lines, i.e.\ $G_{Z_{d,g}}$ has
$d-2g+1$ edges;
\item $Z_{d,g}$ has $d-2g +2 $ points of type $R_3$ and $2g-2$ points
of type $S_4$;
\item for each $i$, $V_i$ is the central plane through a point $p$ of
type either $R_3$ or $S_4$;
\item there exist two $R_3$-points of $Z_{d,g}$ whose
central planes do not meet, unless $g=1$ and $d=5$;
\item $\chi(\Oc_{Z_{d,g}})=1-g$, $p_\omega(Z_{d,g})=0$,
$q(Z_{d,g})=g(Z_{d,g})=g$.
\end{itemize}
}

\vskip 6pt
The base of the induction is the case $g = 1$. In this case, $Z_{d,1}$ is
the surface $X_{d,1}$ considered in Construction \ref{ex:3constr1s}.
Now we assume $g >1$ and we describe the inductive step.

Consider the surface $Z_{d-3, g-1}$, which sits in $\pp^{d-2g}$, which
we suppose to be embedded as a hyperplane in $\Pp^{d-2g+1}$.

If $g =2$ and $d=8$, choose two distinct planes $V_1$ and $V_2$
of $Z_{5,1} = X_{5,1}$, which
do not meet along a line. Otherwise, choose two distinct planes $V_{1}$ and $V_{2}$
of $Z_{d-3,g-1}$ which are central for
two $R_3$-points, say $p_{1}$ and $p_{2}$, and which span a $\Pp^5$.

Choose a line $l_1$ in $V_{1}$ [resp.\ $l_2$ in $V_{2}$]
which is general among those lines passing through $p_{1}$ [resp.\ through
$p_{2}$]. Consider a general $\Pp^4$ in $\Pp^{d-2g +1}$ containing $l_1$ and $l_2$.

One can show that, in this $\pp^4$,
there is a smooth, rational normal cubic scroll $R'$ which contains
$l_1$ and $l_2$ and
such that $R'$ meets transversally $Z_{d-3,g-1}$ along $R'\cap
Z_{d-3,g-1}=l_1\cup l_2$.

In this $\pp^4$, the cubic scroll $R'$ degenerates to a planar Zappatic surface
$X_{3,0}$, consisting of three planes, say $V_{d-2}$, $V_{d-1}$ and $V_{d}$,
such that $l_1\subset V_{d}$ and $l_2\subset V_{d-2}$.

We define $Z_{d,g}= Z_{d-3,g-1}\cup X_{3,0}$.
We added three planes and four double lines; the points $p_1$
and $p_2$ becomes of type $S_4$
for $Z_{d,g}$ and we added three points of type $R_3$ at $V_1\cap
V_{d-1}\cap V_{d}$, at $V_{2}\cap V_{d-2}\cap V_{d-1}$
and at $V_{d-2}\cap V_{d-1}\cap V_{d}$.
It is clear the existence of two $R_3$-points whose central planes do not meet.

Arguing by induction, one may assume that $Z_{d-3,g-1}$ is the flat limit of a smooth,
linearly normal scroll $S$ of degree $d-3$ and genus $g-1$ in
$\Pp^{d-2g}$. If $F_i$, $i=1,2$, is
the ruling of $S$ whose limit is $l_i$ and $R$ is a
smooth, cubic scroll containing $F_1$, $F_2$ as ruling
and whose limit is $R'$,
one can show, by using the same proof of Theorem \ref{thm:deg3},
that the union $S\cup R$ is the flat limit of a family of smooth, linearly
normal scrolls of degree $d$ and genus $g$ in $\Pp^{d-2g+1}$,
which is contained in the same component ${\mathcal H}_{d,g}$ of
Theorem \ref{thm:deg3} (cf.\ Theorem \ref{thm:contoflam} and Remark \ref{rem:new}).
\end{example}

\section{Hilbert schemes of scrolls}\label{S:compo}

In this section we prove that ${\mathcal H}_{d,g}$, as determined
in Theorem \ref{thm:deg3}, is the unique irreducible component of
the Hilbert scheme of scrolls of degree $d$ and genus $g$ in
$\Pp^{d-2g+1}$ whose general point parametrizes a smooth, linearly
normal scroll (cf.\ Theorem \ref{thm:contoflam}). This component
$\HH_{d,g}$ dominates
${\mathcal M}_g$ (cf.\ Remark \ref{rem:hilbconto}).

This, together with Construction
\ref{ex:3constr} and Theorem \ref{thm:deg3},
proves Theorem \ref{thm:intro} in the introduction.

On the other hand, we will also construct
families of scrolls $Y$ of degree $d$ and genus $g$
in $\Pp^r$, with $r > d-2g +1$, with
$h^1(Y, \Oc_Y(1))\neq 0$ (cf.\ Example \ref{ex:contoflam2}). We will also show
that projections of such scrolls may fill up components of the Hilbert scheme,
different from $\HH_{d,g}$, which may even dominate ${\mathcal M}_g$ (cf.\ Example \ref{ex:ciro}).

Let $C$ be a smooth curve of genus $g$ and let $F \stackrel{\rho}{\to} C$ be a {\em geometrically
ruled surface} on $C$,
i.e. $F = \Pp(\Ff)$, for some rank-two vector bundle $\Ff$ on $C$. Furthermore, we assume that
$\Ff$ is very ample, i.e.
$F$ is embedded in $\Pp^r$,
for some $r \geq 3$, via the $\Oc_F(1)$ bundle as a scroll of degree $d= deg(\Ff)$. From now on, $H$ will
denote the hyperplane section of $F$. A general hyperplane section $H$ is isomorphic to $C$,
so that we will set $L_F$ the line bundle on $C \cong H$ which is the restriction of the hyperplane bundle.
We will denote by $R$ a general ruling of $F$, and more precisely by $R_x$ the ruling mapping to the point $x$
in $C$.

Let $Y := C \times \Pp^1$. If $L$ is a line bundle on $C$, we will set
\begin{equation}\label{eq:sqtimes}
\tilde{L}:=  \pi_1^*(L) \otimes \pi_2^*(\Oc_{\Pp^1}(1)),
\end{equation}where $\pi_i$ denotes the projection on the $i^{th}$-factor, $1 \leq i \leq 2$.

\begin{proposition}\label{prop:ghio}
Let $C$ be a smooth curve of genus $g \geq 0$ and let $F:= \Pp(\Ff)$ be a geometrically
ruled surface on $C$. Assume that
$deg(\Ff) = d$.

Then there is a birational map$$\varphi : Y \dasharrow F$$which is the composition
of $d$ elementary transformations at distinct points of a set
$\Gamma:= \{ y_1, \ldots, y_d\} \subset Y$ lying on $d$ distinct rulings of $Y$. Moreover,
\begin{itemize}
\item[(i)] $\varphi^*(\Oc_F(H)) = \tilde{L_F}$; \item[(ii)]
$\varphi^*(|\Oc_F(H)|) = |\tilde{L_F} \otimes \Ii_{\Gamma/Y}|$.
\end{itemize}
\end{proposition}
\begin{proof}
The argument is similar to the one in \cite{Ghio}, Prop. 6.2, and in \cite{MN}. Indeed, let $\Pi$ be a general
linear
subspace of codimension two in $\Pp^r$ which is the base locus of a pencil ${\mathcal P} \cong \Pp^1$ of hyperplanes.
By abusing notation, we will denote by $ {\mathcal P}$ the corresponding pencil of hyperplane sections of $F$.
More specifically, we will denote by $H_t$ the hyperplane section corresponding to the point $t \in \Pp^1$.
Then we denote by $Z := \{ z_1, \ldots, z_d \} = F \cap {\mathcal P}$; note that
$Z$ is formed by distinct points on distinct rulings.

The map $\varphi:  Y \dasharrow F$ is defined by sending the general point $(x,t) \in Y$ to the point $R_x
\cap H_t \in F$. One verifies that $\varphi$ is birational and  that the indeterminacy locus on $F$ is $Z$.
In order to describe
the map $\varphi$ on $Y$, note that each point $z_i$ maps to a point $x_i \in C$ and determines a unique
value $t_i \in \Pp^1$ such that $H_{t_i}$ contains the ruling $R_{x_i}$,  $1 \leq i \leq d$.
The indeterminacy locus of $\varphi$ on $Y$ is $\Gamma:= \{ y_1, \ldots, y_d\}$, where $y_i = (x_i, t_i)$,
$1 \leq i \leq d$.

As shown in \cite{Ghio}, $\varphi$ is the composition of the elementary transformations based at the points of
$\Gamma$. The rest of the assertion immediately follows.
\end{proof}

Let $\Gamma = \{ y_1, \ldots, y_d \} \subset Y$ be a subset formed by $d$ distinct points. We consider the
line bundle on $C$
\begin{equation}\label{eq:lgamma}
L_{\Gamma} :=\Oc_{C}( x_1 + \ldots + x_d),
\end{equation}where $\pi_1(y_i) = x_i$, $ 1 \leq i \leq d$.

\begin{theorem}\label{thm:contoflam}
Let $g \geq 0$ and $d > 2g+3$ be integers. Then there exists a
unique irreducible component $\HH_{d,g}$ of the Hilbert scheme,
parametrizing scrolls of degree $d$ and genus $g$ in
$\Pp^{d-2g+1}$, whose general point represents a smooth scroll $F
\subset \Pp^{d-2g+1}$ which is linearly normal and moreover with
$h^1(F, \Oc_F(1)) = 0$.
\end{theorem}
\begin{proof}
Let $U \subset Hilb^d(Y)$ be the open subset formed by all $\Gamma = \{ y_1, \ldots, y_d \} \subset Y$
containing $d$ points lying on $d$ distinct fibres and imposing $d$ independent conditions on $| \tilde{L_{\Gamma}}|$,
which means$$dim ( | \tilde{L_{\Gamma}} \otimes \Ii_{\Gamma/Y} |) = dim (| \tilde{L_{\Gamma}}|) -d.$$Note that, by
the Kunneth formula,
$h^0( \tilde{L_{\Gamma}}) = 2 h^0(L_{\Gamma}) = 2(d-g+1)$.
Thus, $dim ( | \tilde{L_{\Gamma}} \otimes \Ii_{\Gamma/Y} |) = d-2g+1.$ The linear system
$|\tilde{L_{\Gamma}}|$ determines a rational map$$\varphi : Y \dasharrow \Pp^{d-2g+1}.$$By Proposition \ref{prop:ghio},
every smooth scroll $F$ of degree $d$ and genus $g$ in $\Pp^{d-2g+1}$ is the image of such a map. Therefore,
for general $\Gamma$ in $U$, the map $\varphi$ is birational onto its image $F$, which is a smooth
scroll of degree $d$ and genus $g$ whose Hilbert point $[F]$ belongs to a unique
well-determined component $\HH_{d,g}$ of the Hilbert scheme.

Note that by $(ii)$ of Proposition \ref{prop:ghio}, $h^1 (\Oc_F(1)) = h^1( \tilde{L_{\Gamma}} \otimes \Ii_{\Gamma/Y} )
= 0$; therefore, by the Riemann-Roch Theorem, $h^0(\Oc_F(1)) = d -2 g+2.$
\end{proof}

\begin{remark}\label{rem:new} Observe that
the irreducible component $\HH_{d,g}$ determined in Theorem \ref{thm:deg3} coincides
with the one determined in Theorem \ref{thm:contoflam}. The case $d = 2g+3$ can also be covered
with similar arguments. In that case, we have surfaces in $\Pp^4$ which are no longer smooth,
but they have $2g(g-1)$ double points as dictated by the {\em double point formula}. Nonetheless, the
statement of Theorem \ref{thm:contoflam} still holds by substituting $F$ with its normalization.
\end{remark}

\begin{remark}\label{rem:hilbconto} The dimension count for $\HH_{d,g}$ which has been done in Thereom \ref{thm:deg3}
also stems from the proof of Theorem \ref{thm:contoflam}, which provides a parametric representation of $\HH_{d,g}$.
Indeed, the number of parameters on which the general point of $\HH_{d,g}$ depends, is given by the following
count:
\begin{itemize}
\item $3g-3$ parameters for the class of the curve $C$ in ${\mathcal M}_g$, plus
\item $2d$ parameters for the general point in $U$, plus
\item $(r+1)^2 -1$ parameters for projective transformations in $\Pp^r$, where $r = d-2g+1$, minus
\item $2(r-1) = 2d -4g $ parameters for the choice of a codimension-two subspace $\Pi$ in $\Pp^r$, minus
\item $3$ parameters for projective isomorphisms of the pencil of hyperplanes through $\Pi$ with $\Pp^1$.
\end{itemize}This computation shows that $\HH_{d,g}$ has {\em general moduli}, in the sense that
the base of the general scroll $[F] \in \HH_{d,g}$ is a general point of ${\mathcal M}_g$.

Observe that this can also be viewed as  a consequence of Theorem \ref{thm:deg3} and more specifically
of the fact that $h^1(\Oc_F(1)) = 0$ for $[F]$ a general point of the generically smooth component
$\HH_{d,g}$.

Indeed, if $F \subset \Pp^r$, $r
= d-2g+1$, is a smooth scroll, from the Euler sequence restricted to
$F$,
\[
0 \to \Oc_F \to H^0(\Oc_F(1))^{\vee} \otimes \Oc_F(1) \to
\T_{\Pp^r|F} \to 0,
\]we get that $h^1(\T_{\Pp^r|F}) = 0$. Therefore,
from the normal sequence of $F$ in $\Pp^r$
\[
0 \to \T_F \to \T_{\Pp^r|F} \to \N_{F/\Pp^r}\to 0,
\]we get the surjection
\[
H^0( \N_{F/\Pp^r})\to\!\!\!\to  H^1(\T_F).
\]Since $F$ is a $\Pp^1$-bundle over $C$, from the differential
of the map $F \stackrel{\rho}{\to} C$, we get a surjection
\[
H^1(\T_F)\to\!\!\!\to  H^1(\T_C),
\]hence
\[
H^0( \N_{F/\Pp^r})\to\!\!\!\to  H^1(\T_C).
\]
which shows that $\HH_{d,g}$ dominates ${\mathcal M}_g$.
\end{remark}

Next, we consider the problem of the existence of components of the Hilbert schemes of scrolls
of degree $d$ and genus $g$ in $\Pp^r$, with $r > d-2g+1$. First, it is easy to determine an upper-bound for $r$.
This subject has been deeply studied by C. Segre (cf.\ \cite{Seg} and \cite{Ghio}). For the following lemma, compare
\cite{Seg}, \S\ 14.

\begin{lemma}\label{lem:segre} Let $g \geq 1$ be an integer.
Let $C$ be a smooth curve of genus $g$ and let $F = \Pp(\Ff)$ be a ruled surface on $C$ and $d = deg(\Ff)
\geq 2g+1$. Assume that there exists a smooth curve in $| \Oc_F(1) |$. Then,$$h^0(\Oc_F(1))
\leq d -g + 2.$$The equality holds if and only if $\Ff = \Oc_C \oplus L$,
in which case $\Oc_F(1) $ maps $F$ to a cone over a projectively normal curve of degree $d$ and genus $g$
in $\Pp^{d-g}$.
\end{lemma}
\begin{proof}
The bound on $h^0(\Oc_F(1))$ follows by the Riemann-Roch Theorem on $C$. If the equality holds, then $C$
is linearly normally
embedded as a curve of degree $d$ and genus $g$ in $\Pp^{d-g}$. It is well-known that this curve is projectively normal
(cf.\ \cite{Cast}, \cite{Matt} and \cite{Mumf}).
Therefore $F$ is mapped to a surface $X$ which is projectively normal, since its
general hyperplane section is (cf.\ \cite{Greco}, Theorem 4.27).

On the other hand, $X$ is a scroll of positive genus. Therefore $X$ cannot be smooth, and it has
some isolated singularities. This forces $X$ to be a cone (cf.\ Claim 4.4 in \cite{CLM2}). Hence, the assertion follows.
\end{proof}

\begin{remark} Let $C$ be a smooth curve of genus $g$ and let $F = \Pp(\Ff)$ be a ruled surface on $C$ and $d = deg(\Ff)
\geq 2g+1$. Then
\begin{equation}\label{eq:boundh0}
d - 2g + 2 \leq  h^0( \Oc_F(1)) \leq d -g + 2,
\end{equation}where the lower bound is immediately implied by the Riemann-Roch Theorem whereas the upper bound is given
by the previous lemma. Equivalently,
\begin{equation}\label{eq:boundh1}
0 \leq  h^1( \Oc_F(1)) \leq g,
\end{equation}where the upper-bound is realized by the cones and the lower-bound by the general scrolls in the
component $\HH_{d,g}$ considered above.

Any intermediate value $i$ of $h^1(\Oc_F(1))$, $1 \leq i \leq g$, can be actually realized.
An easy construction is via decomposable bundles as the following example shows.
\end{remark}

\begin{example}\label{ex:contoflam2}
Let $g \geq 3$ and let $d \geq 4g -1$ be integers. Let $i$ be any integer between $1$ and $g$.
Let $C$ be a smooth, projective curve of genus $g$ with a line bundle $L$ such that $|L|$ is base-point-free
and $h^1(L) = i$. Let $D$ be a general divisor of degree $d - deg (L)$. Notice that, since
$deg(L) \leq 2g-2$ and $d \geq 4g -1$, then  $deg(D) \geq 2g+1$ and the linear series
$|D|$ is very ample.

Consider $\Ff = L \oplus \Oc_C(D)$. If $F = \Pp(\Ff)$ then $\Oc_F(1)$ is base-point-free and
$h^1(\Oc_F(1)) = i$.

For large values of $i$, $\Oc_F(1)$ is never very ample. For instance, for $i = g-1$, $C$ is forced to be
hyperelliptic and $L = g^1_2$. Thus, the image of $F$ via $|\Oc_F(1)|$ has a double line.

Similarly, if $i = g-2$, either $C$ is hyperelliptic and $L = 2 g^1_2$, or $C$
is trigonal and $L = g^1_3$ or $g=3$ and $L = \omega_C$. In the former case, the image of $F$ has a double
conic; in the second case, the image of $S$ has a triple line. Only in the third case, the image of $C$ via
$|L|$ is smooth.

The analysis is subtle and we do not dwell here on this.
\end{example}

Now we consider the question of whether there are other components,
different from $\HH_{d,g}$,  of the Hilbert scheme of surfaces in $\Pp^{d-2g+1}$ whose general point corresponds to
a smooth scroll of degree $d$ and genus $g$. The answer to this question is affirmative; in fact one can
construct such components even with general moduli. In the next example, we show one possible
construction of a component with general moduli. The reader may easily generate other
similar constructions.

\begin{example}\label{ex:ciro}
Let $C$ be a curve with general moduli of genus $g = 4l + \epsilon$, where
$ 0 \leq \epsilon \leq 3$. Let $L$ be a very-ample, special line bundle of
degree $m:= 3+g-l$ with $h^0(L) = 4$. Note that such a $L$ varies in a family
of dimension $\rho := \rho (g, 3,m) = \epsilon$.

Let $d$ be an integer with either $d \geq 2g + 10$, if $ \epsilon = 0, 1$, or $d \geq 2g + 11$, if
$\epsilon = 2, 3$. Set $r = d-2g+1$.

Let $N$ be a general line bundle on $C$ of degree $d-m$. Note that $d-m > g+7+l$. Hence $N$ is very ample (cf.\ e.g.
\cite{AC}) and $h^0(N) = d-m-g +1$.

Set $\Ff = L \oplus N$ and  $X = \Pp(\Ff)$.
Then $R+1 := h^0(\Oc_X(1)) = h^0(L) + h^0(M) = r + 1 + l$.

Since $\Oc_X(1)$ is very ample, $X$ is linearly normal embedded in $\Pp^R$ as a smooth scroll
of degree $d$ and genus $g$, which can be generically projected to $\Pp^r$ to a smooth scroll
$X'$ with the same degree and genus, which belongs a certain component $\HH$ of the Hilbert
scheme. As in the proof of Theorem \ref{thm:deg3}, the general member of $\HH$ is a scroll
of the same degree and genus.

The dimension of $\HH$ can be easily bounded from below by the sum of the following quantities:
\begin{itemize}
\item $3g -3$, which are the parameters on which $C$ depends,
\item $g$, which are the parameters on which $N$ depends,
\item $\epsilon$, which are the parameters on which $L$ depends,
\item $(r+1)l = dim (\G (r, R)) $, which are the parameters for the projections,
\item $(r+1)^2 -1 = dim (PGL(r+1, \CC)).$
\end{itemize}

The hypothesis on $d$ implies that  $dim(\HH) \geq dim (\HH_{d,g})$, which shows that $\HH$ is different from
$\HH_{d,g}$.
\end{example}

\begin{remark} The question of understanding how many components of the Hilbert scheme of
scrolls there are, and the corresponding image to the moduli space of curves of genus $g$,
is an intriguing one. The previous example suggests that the a complete answer could be rather complicated.
It also leaves open the question whether $\HH_{d,g}$ is the only component with general moduli
for $ 2g + 4 \leq d \leq 2g + 10$.
\end{remark}

\section{Comments on Zappa's original approach}\label{S:Zappa}

In \cite{Zap1}, Zappa stated a result about
embedded degenerations of scrolls of
sectional genus $g \geq 2$ to unions of planes. His result, in our terminology, reads
as Theorem \ref{thm:zappa} in the introduction.

Zappa's arguments rely on a rather intricate analysis of algebro-geometric and topological type
of degenerations of hyperplane sections of the scroll and, accordingly, of the branch curve of a general
projection of the scroll to a plane.

We have not been able to check all the details of this very clever argument. This is
one of the reason why we preferred to solve the problem in a different way, which is the one
we exposed in the previous sections. Our approach has the advantage of proving a result
in the style of Zappa, but with better hypotheses about the degree of the scrolls.

However, the idea which Zappa exploits, of degenerating the branch curve of a general projection to a plane,
is a classical one which goes back to Enriques, Chisini, etc, and certainly deserves attention.
We hope to come back to these ideas in the future.

In reading Zappa's paper \cite{Zap1}, our attention has been attracted also by another
ingredient he uses which looks interesting on its own. It gives extendability
conditions for a curve on a scroll which is not a cone.
We finish this paper by briefly
reporting on this. At the the end of the section we briefly summarize Zappa' s argument
for the degenerations of the scroll.

Let $F\subset \pp^3$ be a scroll, which is not a cone over a plane curve.
We do not assume $F$ to be smooth. Equivalently, we can look at $F$ as a curve
$\CC$ in the Grassmannian
$\G(1,3)$ of lines in $\pp^3$, which
is isomorphic to the Klein hyperquadric in $\pp^5$ via the Pl\"ucker embedding.

Let $\Pi$ be a general plane and let
$\Gamma:=F\cap \Pi$. Consider
$$\nu:C\to\Gamma$$the normalization map. Then, there is a commutative
diagram
\[
\xymatrix@C-10pt{%
C \ar[rr]^-{\Phi} \ar[dr]!UL_-{\nu}
 & & \CC \subset \G(1,3) \subset \pp^5 \ar[]!DL;[dl]!UL+<1.2em,0ex>^-{\pi}  \\
& \Gamma \subset \Pi
}
\]where $\Phi$ maps a general point $x \in C$ to the unique line of $F$ passing through $\nu(x)$,
and $\pi$ maps each point $l \in \CC$, corresponding to a ruling $L$ of $F$, to the point
$L\cap \Gamma$.

Zappa proves the following nice lemma:

\begin{lemma}\label{lem:omog}(cf.\ \S 1 in \cite{Zap1}) In the above setting:
$$\nu^*(\Oc_{\Gamma}(1) ) \cong \Phi^*(\Oc_{\CC}(1)).$$More specifically,
$\pi$ is the projection of $\CC$ from the plane $\Pi^* \subset \G(1,3)$,
filled up by all lines of $\Pi$.
\end{lemma}
\begin{proof}The assertion follows
from the fact that, if $r$ is a line in $\Pi$, then $\pi^*(r) $ is the section of the tangent hyperplane to
$\G(1,3)$ at the point of $\Pi^*$ corresponding to $r$. Such a hyperplane contains $\Pi^*$, and conversely any
hyperplane containing $\Pi^*$ is of this type.
\end{proof}

Zappa notes that an interesting converse of the previous lemma holds.

\begin{prop}\label{p:1} (cf. \S\ 2 in \cite{Zap1})
An irreducible plane curve $\Gamma$ is a section of a scroll $F \subset \Pp^3$ of
degree $d$ if and only if $\Gamma$ is the projection of
a curve $\CC$ of degree $d$, lying on a smooth quadric
$Q\subset\pp^5$, and the center of the projection is a plane
contained in $Q$.
\end{prop}

\begin{proof} One implication is Lemma \ref{lem:omog}. Let us prove the other
implication.

Suppose that  $\Gamma$ is the projection
of $\CC\subset Q\subset \pp^5$ from a plane $\bar\Pi\subset Q$.
Since all smooth quadrics in $\pp^5$ are projectively equivalent,
we may assume that $Q$ is the Klein hyperquadric. The assertion follows by reversing
the argument of the proof of Lemma \ref{lem:omog}.
\end{proof}

Proposition \ref{p:1} can be extended in the following way.
Let $\Gamma$ be a plane curve of degree $d$ and geometric genus $g$, such
that $d\ge g+6$. Set $i=h^1(C,\nu^*({\Gamma}(1)))$.
Then, one has the birational morphism
\begin{equation}\label{eq:emb}
C \xrightarrow{\,|\nu^*(\Oc_{\Gamma}(1))|\,} \bar\CC \subset \pp^r,
\end{equation}where $ r = d-g + i > 5$ and the following linear projection:
\[\bar\CC \stackrel{\bar \pi}{\to} \Gamma \subset\pp^2.
\]

\begin{prop}\label{p:2}(cf.\ \S\ 3 in \cite{Zap1}) In the above setting,
$\Gamma$  is a plane section of
a scroll $F$ in $\pp^3$, which is not a cone, if and only if
$\bar\CC$ lies on a quadric of rank $6$ in $\pp^{r}$ which
contains the center of the projection $\bar \pi$.
\end{prop}

\begin{proof} This is an immediate consequence of Proposition \ref{p:1} and can be left to the reader.
\end{proof}

Zappa uses Proposition \ref{p:2} to prove that any plane curve of degree $d > > g$ is the plane section
of a scroll $F$ which is not a cone. The next
proposition is essentially Zappa's result in \S\ 7 of
\cite{Zap1}, with an improvement on the bound on $d$: Zappa's bound is
$d \geq 3g+2$.

\begin{lemma}\label{p:3} Let $g \geq 0$ and let $d \geq 2g+2$ be integers.
Let $\bar\CC$ be an irreducible, smooth curve of degree $d$ and genus $g$
in $\pp^r$, $r = d-g$. Then there exists a
quadric of $\pp^r$, of rank at most 6, which contains
$\bar\CC$ and a general $\pp^{r-3}$.
\end{lemma}
\begin{proof}
Note that a quadric $Q$ of $\pp^r$ contains a $\pp^{r-3}$ if and only if
$Q$ has rank at most 6.

Consider the short exact sequence
\[
0\to \II_{\bar\CC/\pp^r}(2) \to \OO_{\pp^r}(2) \to \OO_{\bar\CC}(2) \to 0.
\]
Since $d \geq 2g+2$, one has  $h^0(\OO_{\bar\CC}(2))=2d-g+1$ and $\bar{\CC}$ is
projectively normal (cf.\ \cite{Cast}, \cite{Matt}, \cite{Mumf}).
Thus
\begin{equation}\label{eq:ciro}
h^0(\II_{\bar\CC/\pp^r}(2)) = \binom{r+2}{2} - (2d-g+1).
\end{equation}Let $\Sigma$ be a general $\Pp^{r-3}$ in $\Pp^r$.
Then, from \eqref{eq:ciro}, one has
\[
h^0(\II_{\bar{\CC} \cup \Sigma /\pp^r}(2)) \geq  \binom{r+2}{2} - \binom{r-1}{2} -
(2d-g+1) = d - 2g -1 > 0.
\]
\end{proof}

We need the following lemma:

\begin{lem}\label{l:1}
Let $\bar\CC\subset\pp^r$ be as in Proposition \ref{p:3} and assume that, if $g = 0$,
$d \geq 3$. Let $\Sigma$ be a $\Pp^{r-3}$.
The general quadric in the linear system $|\II_{\bar{\CC} \cup \Sigma /\pp^r}(2)|$
has rank $k> 3$.
\end{lem}

\begin{proof}
Suppose by contradiction that all quadrics containing $\bar\CC$
and $\Sigma$ have rank 3.
Let us define
\[
R_3(\bar\CC):=\{ Q\in \pp(H^0(\II_{\bar C/\pp^r}(2))) \, | \,
\rk(Q) \le 3 \}.
\]By an easy count of parameters our assumption implies that:
\[
\dim R_3(\bar\CC)\ge 3d-4g-7.
\]Next, we will show that this inequality is not possible.

In order to do that, we apply results from
\cite{Zam}. Zamora proves in \cite{Zam}, cf.\ Lemma 1.2, that there is a one-to-one correspondence
between quadrics $Q \in R_3(\bar \CC)$ and pairs
$(g^1_a,g^1_b)$ of linear series on $\bar \CC$, with $a\le b$, such that:
\begin{enumerate}[(i)]
\item $a+b=\deg \bar\CC=d$,
\item $|g^1_a+g^1_b|=|\OO_{\bar\CC}(1)|$,
\item $g^1_a+B_b=g^1_b+B_a$, where $B_a$ ($B_b$, resp.) is the base locus
of the $g^1_a$ ($g^1_b$, respectively).
\end{enumerate}

Let $Q$ be the general member of an irreducible component $W$ of maximal dimension of
$R_3(\bar \CC)$ and let $(g^1_a, g^1_b)$ be the corresponding pair of linear series on
$\bar \CC$.

Zamora's result implies that there is a base-point-free linear series $g^1_h$ on $\bar \CC$
such that\[ g^1_a  = g^1_h + B_a, \; g^1_b = g^1_h + B_b,\]so that
$$|\Oc_{\bar \CC}(1) | = | 2 g^1_h + B_a + B_b |.$$Note that, once the divisor
$B_a + B_b$ has been fixed, the line bundle $L$ corresponding to $g^1_h$ belongs to a
zero-dimensional set in $Pic^h(\bar \CC)$. Set $\delta = deg (B_a + B_b)$, so that
$d = 2h + \delta$.

Suppose now that $L$ is non-special. Then,
$$3d - 4g -7 \leq dim (W) \leq \delta + 2 (h - g - 1) = d - 2g -2,$$which gives a contradiction.

Now assume that $L$ is special, so that $|L| = g^r_h$, with $2 r \leq h$. In this case
$$3d - 4g -7 \leq dim (W) \leq \delta + 2 (r- 1) \leq  \delta + h - 2,$$which
leads to a contradiction.
\end{proof}

As a consequence of the previous lemma, we have:

\begin{thm}\label{t:1}
Let $\Gamma$ be an irreducible, plane curve of degree $d$ and geometric genus $g
\geq 0$. If $d\ge {\rm max}\ \{g+5, 2g+2 \}$, then $\Gamma$ is a plane section of a
scroll in $\pp^3$, which is not a cone.
\end{thm}

\begin{proof}
Let $\bar\CC\subset\pp^r$ be the curve corresponding to $\Gamma$
in $\pp^2$. Then $\Gamma$ is the projection of $\bar\CC$ from $\Sigma = \Pp^{r-3}$
disjoint from $\bar\CC$. By Lemma \ref{p:3}, there is a quadric $Q$ containing $\bar\CC
\cup \Sigma$. If $rank (Q) := k$ is $6$, we finished by Proposition \ref{p:2}.
By Lemma \ref{l:1}, we know that $k \geq 4$.

If $k = 5$,  then the vertex $V$ of $Q$ is a $\Pp^{r-5}$. By projecting from $V$, $Q$ maps
to a smooth quadric $Q'$ in $\Pp^4$ containing $\CC'$, the projection of $\bar \CC$, and $\Sigma'$,
the projection of $\Sigma$; the line $\Sigma'$ is skew with respect to $\CC'$. Of course
$\Gamma$ is the projection of $\CC'$ from $\Sigma'$. Let us embedd $\Pp^4$ in $\Pp^5$ as a
hyperplane. We can certainly find a smooth quadric $\bar Q$ in $\Pp^5$ containing $Q'$ and
containing a plane $\Pi$ intersecting the $\Pp^4$ in $\Sigma'$. The curve $\Gamma$ is now the
projection of $\CC'$ from $\Pi$. The assertion follows from Proposition \ref{p:1}.

If $k = 4$, then the vertex $V$ of $Q$ is a $\Pp^{r-4}$. Suppose first that $\Sigma$ contains $V$; then
by projecting from $V$ to $\Pp^3$, the quadric $Q$ maps to a smooth quadric $Q'$, containing
$\CC'$, the image of $\bar \CC$, and the point $p \in Q'$, the image of $\Sigma$, which
does not sit on $\CC'$. The curve $\Gamma$ is the projection of $\CC'$ from $p$. At this point, we can finish as
in the previous case, by embedding the $\Pp^3$ in $\Pp^5$ and finding a smooth quadric $\bar Q$ in
$\Pp^5$ containing $Q'$ and the plane $\Pi$ intersecting the $\Pp^3$ in $p$.

If $\Sigma$ does not contain $V$, it intersects $V$ in $W \cong \Pp^{r-5}$. By projecting from
$W$ to $\Pp^4$ we get a situation similar to the case $k=5$. The only difference is that
$Q'$ is now singular at a point $p$, however $\Sigma'$, the projection of $\Sigma$, does
not contain $p$. So we can conclude exactly as in the case $k=5$.
\end{proof}

\begin{remark}\label{c:1}
We add a little remark to Theorem \ref{t:1}.
Let $\Gamma$ be a plane curve which is a plane section of a scroll $F \subset \Pp^3$, which is
not a cone. So if one applies Theorem \ref{t:1}, the scroll which extends $\Gamma$ is certainly
not developable.
\end{remark}

As Zappa does in \cite{Zap1}, one can get an
interesting consequence of Theorem \ref{t:1} by applying duality.
Recall that the {\em class} of an irreducible plane curve is the degree of the
dual curve.

\begin{corollary}\label{c:2} An irreducible, plane curve of class $d$ and geometric genus
$g$, such that $d \geq  {\rm max}\ \{g+5, 2g+2 \}$, is the branch curve of a projection of a scroll
in $\Pp^3$ of degree $d$ and genus $g$, which is not a cone.
\end{corollary}
\begin{proof}Let $D \subset \Pp^2$ be an irreducible plane curve of class $d$. Let
$\Gamma \subset (\Pp^2)^*$ be the dual curve. By Theorem \ref{t:1}, $\Gamma$ is the plane section
of a scroll $\Phi$ which is not a cone. By standard properties of duality, $D$ is the branch curve
of the projection of $F = \Phi^*$ from the point corresponding to the plane in which
$\Gamma$ sits.
\end{proof}

The argument of Zappa to prove the degeneration of a scroll to a union of planes
runs as follows. Zappa considers the scroll $F$ whose hyperplane section
$\Gamma$ is a general member
of the Severi variety $V_{d,g}$ of plane curves of degree $d$ and geometric
genus $g$. Then he lets $\Gamma$ degenerate to a general union of $d$ lines. From a
complicated analysis involving the degeneration of $\Gamma$ and the degeneration of its dual
curve, which is the branch curve of the projection of the dual of the surface on the plane
(see Corollary \ref{c:2}), Zappa deduces that in this degeneration of $\Gamma$, $F$
degenerates to a union of planes. Moreover, he controls the degeneration of the
linearly normal model of $F$ deducing that it also degenerates to a union of planes with only
points of type $R_3$ and $S_4$.



\end{document}